\documentclass[11pt]{article}

\usepackage[T1]{fontenc}
\usepackage{amsmath,amssymb,amsthm}
\usepackage{mathtools}
\usepackage[margin=1.1in]{geometry}
\usepackage{tikz}
\usepackage{float}
\usepackage{placeins}

\usepackage[hidelinks]{hyperref}

\theoremstyle{definition}
\newtheorem{definition}{Definition}
\theoremstyle{plain}
\newtheorem{lemma}[definition]{Lemma}

\newtheorem{theorem}[definition]{Theorem}
\newtheorem{corollary}[definition]{Corollary}
\theoremstyle{remark}
\newtheorem{remark}[definition]{Remark}

\newcommand{\lt}{\lhd_{t}}
\newcommand{\rt}{\blacktriangleleft_{t}}
\newcommand{\lo}{\lhd_{1}}
\newcommand{\ro}{\blacktriangleleft_{1}}
\newcommand{\ordl}[1]{\lhd_{#1}}
\newcommand{\ordr}[1]{\blacktriangleleft_{#1}}

\newcommand{\lmap}[1]{\ell_{#1}}
\newcommand{\rmap}[1]{r_{#1}}

\newcommand{\Nt}{N_{t}}
\newcommand{\No}{N_{1}}
\newcommand{\Nc}[1]{N^{\circ}_{#1}}

\newcommand{\negt}{\neg_{t}}
\newcommand{\nego}{\neg_{1}}
\newcommand{\simt}{\mathord{\sim}_{t}}
\newcommand{\simo}{\mathord{\sim}_{1}}

\newcommand{\forc}[1]{\models_{#1}}
\newcommand{\nforc}[1]{\not\models_{#1}}

\newcommand{\Pset}[2]{P^{#1}_{#2}}

\newcommand{\entt}{\models^{t}_{oe}}
\newcommand{\entone}{\models^{1}_{oe}}
\newcommand{\entn}{\models^{ndst}_{oe}}
\newcommand{\dert}{\vdash_{t}}
\newcommand{\derone}{\vdash_{1}}
\newcommand{\dern}{\vdash_{n}}

\newcommand{\For}{\mathit{For}_{oe}}
\newcommand{\Var}{\mathit{Var}}
\newcommand{\Lang}{\mathcal{L}_{oe}}
\newcommand{\NDst}{\mathit{NDst}_{oe}}
\newcommand{\NDstt}{\mathit{NDst}^{t}_{oe}}
\newcommand{\NDstone}{\mathit{NDst}^{1}_{oe}}

\newcommand{\Frm}{\mathfrak{F}}
\newcommand{\Str}{\mathfrak{X}}
\newcommand{\Latt}{\mathfrak{S}}
\newcommand{\Lstab}{\mathcal{L}}
\newcommand{\Th}[2]{T^{#1}_{#2}}
\newcommand{\Cth}[2]{T^{\circ #1}_{#2}}
\newcommand{\Fc}{\Frm^{c}}

\newcommand{\Four}{\mathbf{4}_{oe}}
\newcommand{\Nine}{\mathbf{9}_{u}}
\newcommand{\LatFour}{\mathfrak{L}_{4}}
\newcommand{\LatNine}{\mathfrak{L}_{9}}
\newcommand{\FDEoe}{\mathit{FDE}_{oe}}
\newcommand{\FDEoeu}{\mathit{FDE}^{u}_{oe}}
\newcommand{\entu}{\models^{u}_{oe}}
\newcommand{\entfde}{\models_{oe}}
\newcommand{\leqt}{\leqslant_{t}}
\newcommand{\leqo}{\leqslant_{1}}

\newcommand{\entfour}[1]{\models^{\Four}_{#1}}
\newcommand{\defiff}{\leftrightharpoons}

\title{Relational Semantics and Completeness for Onto-Epistemic Logic}
\author{O.~Grigoriev\\[0.5ex]
  \normalsize Department of Logic, Faculty of Philosophy\\
  \normalsize Lomonosov Moscow State University, Moscow, Russia}
\date{}

\begin{document}

\maketitle

\begin{abstract}
This paper investigates a relational semantics for onto-epistemic logic
with two semi-negations. The frames carry two pairs of quasi-orders and
two commuting involutions. Some frames have nondistributive lattices of
stable sets. The four- and nine-valued entailments are recovered as joint
entailment on two particular finite frames. Sound and complete binary
consequence calculi are given for both component entailments and for joint
entailment. The joint Lindenbaum algebra embeds
subdirectly into the product of its two component quotients. This
embedding shows that joint derivability is the intersection of the two
component derivability relations. The finite model property is proved
using Dedekind--MacNeille completion and a filter--ideal representation
of the resulting lattice. All three entailments are decidable by
interleaving the search for a derivation with the search for a finite
countermodel.
\end{abstract}

\section{Introduction}\label{sec:introduction}

Belnap's four epistemic values record what a reasoner has been told
about a sentence. Its ontological truth or falsity depends on how
things are, independently of the reasoner's information. Belnap
combines the four epistemic values with two ontological values to
obtain eight pairs~\cite[p.~47]{BelnapComputer}.

In the four-valued semantics of Zaitsev and
Shramko~\cite[Sections~3--5]{bib96}, the ontological and epistemic
coordinates range over $\{f,t\}$ and $\{0,1\}$, respectively. For example,
$\langle t,0\rangle$ represents a sentence that is
ontologically true and epistemically false. Ordering the four pairs
coordinatewise gives a lattice isomorphic to the truth lattice for
first-degree entailment (FDE)~\cite{Dunn1976}. The resulting lattice-order
entailment agrees with FDE on the conjunction--disjunction fragment.
Each of the semi-negations $\negt$ and $\nego$ reverses the order in
its own coordinate and fixes the other. Petrukhin~\cite{Petrukhin2021}
uses the same operations in a four-valued logic of an agent's
acceptance and rejection.

The De Morgan negation of FDE reverses the lattice order. A single
semi-negation preserves strict increases in its unchanged coordinate,
which gives counterexamples to contraposition for lattice-order
entailment. The composition of the two semi-negations is Boolean
complementation on the four-element lattice, so $\negt A$ and
$\nego A$ take complementary values. Their conjunction takes the
least value and their disjunction the greatest. Under lattice-order
entailment, this gives analogues of explosion and excluded middle
\[
\negt A\wedge\nego A\entfde B,
\qquad
B\entfde\negt A\vee\nego A.
\]

The nine-valued semantics places $u_t$ between $f$ and $t$, and $u_1$
between $0$ and $1$. In the product of these chains, either component
or both can be indeterminate. The composition remains an order-reversing
involution, and both semi-negations fix $\langle u_t,u_1\rangle$.
Take $A$ and $B$ to be distinct variables and assign this pair to $A$.
The least value for $B$ refutes the first schema, and the greatest
refutes the second. Grigoriev~\cite[Sections~1--4 and~6]{bib51} presents
the nine-valued construction and outlines a relational generalisation
under additional frame conditions.

The four- and nine-valued lattices are distributive. Bimb\'o and
Dunn~\cite{BimboDunn2001} give a relational semantics and a complete
calculus for conjunction and disjunction over general lattices.
Allwein and Dunn~\cite[Section~6]{AllweinDunn} base their Kripke
semantics for linear logic without exponentials on Urquhart's
representation of bounded lattices~\cite{Urquhart}, which applies to
nondistributive lattices.
In the canonical representation, a point is a maximal disjoint
filter--ideal pair. A formula is true at that point if its value
belongs to the filter, false if it belongs to the ideal, and
indeterminate if it belongs to neither. Or\l{}owska and
Rewitzky~\cite{OrlowskaRewitzky} develop this three-valued approach
for lattice-based logics.

Dzik, Or\l{}owska and van Alten~\cite{DzikOrlowskaVanAlten} use
Urquhart's frames to represent several kinds of lattice negation.
The discrete representation of D\"untsch, Kwuida and
Or\l{}owska~\cite{DuntschKwuidaOrlowska2017} treats weak complementation
and its dual on the same lattice. Almeida~\cite{Almeida2009} uses
canonical extensions to represent lattice negations.
Hartonas~\cite{Hartonas2025} treats antitone involutions within a
two-sorted duality for lattice expansions.

On a common set of points, we use a pair of quasi-orders for each
component and two commuting involutions to interpret the semi-negations.
Truth, falsity, and indeterminacy are determined separately in the two
components. The component entailments preserve ontological and
epistemic truth, respectively, and joint entailment requires both forms
of preservation.

For the four-valued semantics, Wintein and
Muskens~\cite{WinteinMuskens} give sound and complete tableaux for both
components. Closure of both tableaux characterises joint entailment.

Fix a component and erase occurrences of the semi-negation that acts
as the identity there. The lattice of stable sets in that component is
bounded and carries an antitone involution. Such lattices give the
algebraic semantics of paraconsistent quantum logic
(PQL)~\cite[Section~13]{DallaChiaraGiuntini2002}. Thus if the translated
entailment holds in PQL, the original entailment holds in the chosen
component. Conversely, if the translated entailment fails in PQL,
its finite model property gives a countervaluation in a finite
involution lattice. By Lemma~\ref{lem:finite-frame-realisation}, this
lattice embeds into the lattice of stable sets of a finite frame,
giving a countermodel in the chosen component.
Kamide~\cite{Kamide2021} gives a complete calculus for PQL with one
formula on each side of a sequent.

Aguilera and Massas~\cite[Theorems~4.1 and~5.1]{AguileraMassasEx}
axiomatise the intersection of orthologic and the implication-free
fragment of intuitionistic logic by embedding Ex-lattices into products
of ortholattices and Heyting lattices. They use an analogous embedding
for the intersection of orthomodular and intuitionistic logic in a
language with implication~\cite[Theorems~3.1 and~4.2]{AguileraMassasIEx}.
Finite basis results for varietal joins imply that the intersection of
two finitely axiomatised extensions of $\mathrm{FL}_e^+$ is finitely
axiomatisable~\cite[Corollary~4.6]{GalatosJoins}. The same holds over the
disjunction-free fragment of $\mathrm{FL}_e^+$%
~\cite[Corollary~26]{GalatosOlsonRaftery}.

We recover the four- and nine-valued entailments on two particular
finite frames. For the full class of frames, we give separate binary
consequence calculi for all three entailments.
The maps to the two component quotients of the joint Lindenbaum
algebra together reflect its order. We use this fact to prove that
joint derivability is the intersection of the two component derivability
relations. Component completeness then yields completeness of the joint
calculus. All three entailments have the finite model property on
relational frames and are decidable.

\section{Relational semantics}\label{sec:relational}

\subsection{Finite-valued onto-epistemic semantics}\label{sec:finite}

This subsection reviews the finite-valued semantics of $\FDEoe$
\cite{bib96} and $\FDEoeu$ \cite{bib51}. The motivation for distinguishing
ontological and epistemic truth is discussed in \cite{bib96}; here we describe
the two-component value structures, their operations and the associated
consequence relations.

The language $\Lang$ consists of a countable set $\Var$ of propositional
variables, the binary connectives $\wedge$ and $\vee$, the two unary
connectives $\negt$ and $\nego$, called the \emph{ontological} and the
\emph{epistemic} negation respectively, and the auxiliary symbols, a left and
a right parenthesis. The notion of a formula is defined in the standard way;
the set of formulae of $\Lang$ is denoted by $\For$.

A truth value of a formula is a pair whose first coordinate is drawn from an
ontological and whose second from an epistemic base. The simplest choice of
bases, $\{f,t\}$ (ontological) and $\{0,1\}$ (epistemic) ordered by $f\leqt t$ and $0\leqo 1$, gives the
semantics of $\FDEoe$.

\begin{definition}\label{def:four}
Let $\Four=\{f,t\}\times\{0,1\}$. For
$\langle s,u\rangle,\langle v,w\rangle\in\Four$ put
\[
\langle s,u\rangle\leqslant\langle v,w\rangle
\enskip\defiff\enskip
s\leqt v\ \text{ and }\ u\leqo w .
\]
\end{definition}

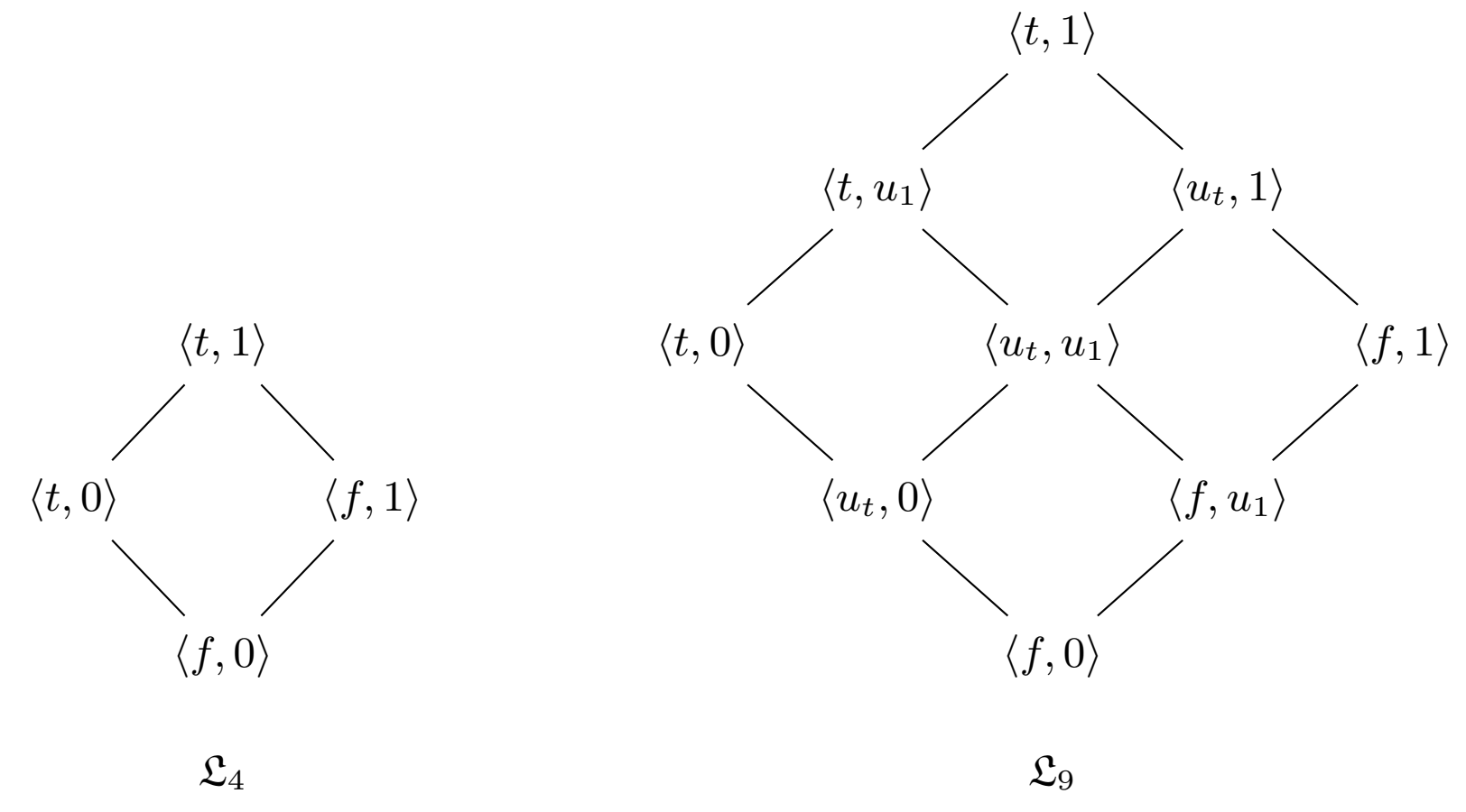
\begin{figure}[ht]
\centering
\begin{tikzpicture}[scale=0.95, every node/.style={font=\footnotesize}]
\node (f0) at (0,0) {$\langle f,0\rangle$};
\node (t0) at (-1.15,1.2) {$\langle t,0\rangle$};
\node (f1) at (1.15,1.2) {$\langle f,1\rangle$};
\node (t1) at (0,2.4) {$\langle t,1\rangle$};
\draw (f0)--(t0) (f0)--(f1) (t0)--(t1) (f1)--(t1);
\node at (0,-0.9) {$\LatFour$};
\begin{scope}[xshift=6.4cm]
\node (a00) at (0,0) {$\langle f,0\rangle$};
\node (a10) at (-1.35,1.2) {$\langle u_t,0\rangle$};
\node (a01) at (1.35,1.2) {$\langle f,u_1\rangle$};
\node (a20) at (-2.7,2.4) {$\langle t,0\rangle$};
\node (a11) at (0,2.4) {$\langle u_t,u_1\rangle$};
\node (a02) at (2.7,2.4) {$\langle f,1\rangle$};
\node (a21) at (-1.35,3.6) {$\langle t,u_1\rangle$};
\node (a12) at (1.35,3.6) {$\langle u_t,1\rangle$};
\node (a22) at (0,4.8) {$\langle t,1\rangle$};
\draw (a00)--(a10) (a00)--(a01);
\draw (a10)--(a20) (a10)--(a11) (a01)--(a11) (a01)--(a02);
\draw (a20)--(a21) (a11)--(a21) (a11)--(a12) (a02)--(a12);
\draw (a21)--(a22) (a12)--(a22);
\node at (0,-0.9) {$\LatNine$};
\end{scope}
\end{tikzpicture}
\caption{The lattices $\LatFour$ and $\LatNine$. In each, the ontological
coordinate increases along the upper-left direction and the epistemic
coordinate along the upper-right one.}\label{fig:lattices}
\end{figure}

The meet and join of each base are written $\sqcap_t,\sqcup_t$ and
$\sqcap_1,\sqcup_1$, and the operation interchanging its two elements is
written $-_t$ and $-_1$ respectively. The operations on $\Four$ are taken
coordinatewise:
\begin{equation}\label{eq:fourops}
\begin{aligned}
\langle s,u\rangle\wedge\langle v,w\rangle
  &=\langle s\sqcap_t v,\ u\sqcap_1 w\rangle,
&\simt\langle s,u\rangle&=\langle {-_t}s,\ u\rangle,\\
\langle s,u\rangle\vee\langle v,w\rangle
  &=\langle s\sqcup_t v,\ u\sqcup_1 w\rangle,
&\simo\langle s,u\rangle&=\langle s,\ {-_1}u\rangle .
\end{aligned}
\end{equation}
Definition~\ref{def:four} makes $\Four$ a lattice, denoted $\LatFour$, with
greatest element $\langle t,1\rangle$, least element $\langle f,0\rangle$ and
two incomparable elements $\langle t,0\rangle$ and $\langle f,1\rangle$; it is
displayed on the left of Figure~\ref{fig:lattices}.

\begin{definition}\label{def:fourval}
A \emph{valuation} of $\Lang$ is a map $v\colon\Var\to\Four$. Its
\emph{extension} is the map $v\colon\For\to\Four$ agreeing with $v$ on $\Var$
and satisfying, for all $A,B\in\For$,
\[
v(A\wedge B)=v(A)\wedge v(B),\quad
v(A\vee B)=v(A)\vee v(B),\quad
v(\negt A)=\simt v(A),\quad
v(\nego A)=\simo v(A),
\]
the operations on the right being those of \eqref{eq:fourops}. Extensions are
denoted by the same letter as the valuations that induce them.
\end{definition}

Several consequence relations may be defined over $\LatFour$, by preservation
of one or of both kinds of truth, or by backward preservation of one of the
kinds of falsity, as is familiar from $\mathrm{FDE}$. We first consider the
relation determined by the lattice order. Its relational generalisation and
the two component relations are introduced in Section~\ref{sec:entailment}.

\begin{definition}\label{def:fdeent}
For $A,B\in\For$,
\[
A\entfde B
\enskip\defiff\enskip
v(A)\leqslant v(B)\ \text{ for every valuation } v .
\]
\end{definition}

Write $t\in v(A)$ to mean that the pair assigned to $A$ by $v$ has $t$ as its
first coordinate, and similarly for $f$, $1$ and $0$. The component
entailment relations are defined by
\begin{equation}\label{eq:fourcomponents}
\begin{aligned}
A\entfour{t}B&\defiff \forall v\,\bigl(t\in v(A)\Rightarrow t\in v(B)\bigr),\\
A\entfour{f}B&\defiff \forall v\,\bigl(f\in v(B)\Rightarrow f\in v(A)\bigr),\\
A\entfour{1}B&\defiff \forall v\,\bigl(1\in v(A)\Rightarrow 1\in v(B)\bigr),\\
A\entfour{0}B&\defiff \forall v\,\bigl(0\in v(B)\Rightarrow 0\in v(A)\bigr).
\end{aligned}
\end{equation}
The superscript marks the four-valued space and is used in this subsection
only; the analogous relations of the relational semantics are introduced in
Section~\ref{sec:entailment} under different notation.

Since each base is two-element, membership of one coordinate excludes
membership of its opposite: $t\in v(A)$ if and only if $f\notin v(A)$, and
$1\in v(A)$ if and only if $0\notin v(A)$. Hence $\entfour{t}$ coincides with
$\entfour{f}$, and $\entfour{1}$ with $\entfour{0}$, so that
\eqref{eq:fourcomponents} yields two relations rather than four. The following
lemma, immediate from Definition~\ref{def:fourval}, states their behaviour
under the connectives.

\begin{lemma}\label{lem:fourcomp}
For all $A,B\in\For$ and every valuation $v$,
\begin{enumerate}
\item[(1)] $t\in v(\negt A)\iff f\in v(A)$ and $f\in v(\negt A)\iff t\in v(A)$;
\item[(2)] $1\in v(\nego A)\iff 0\in v(A)$ and $0\in v(\nego A)\iff 1\in v(A)$;
\item[(3)] $t\in v(A)\iff t\in v(\nego A)$ and $1\in v(A)\iff 1\in v(\negt A)$;
\item[(4)] $t\in v(A\wedge B)\iff t\in v(A)$ and $t\in v(B)$;
\item[(5)] $t\in v(A\vee B)\iff t\in v(A)$ or $t\in v(B)$,
\end{enumerate}
together with the statements obtained from \emph{(4)} and \emph{(5)} by
interchanging $t$ with $f$ and $\wedge$ with $\vee$, and with the statements
obtained from \emph{(1)}--\emph{(5)} by interchanging $t,f,\negt$ with
$1,0,\nego$.
\end{lemma}

Item~(3) is the characteristic feature of the construction: each negation
leaves the coordinate of the other kind untouched, so that the two components
of a value evolve independently under the connectives. The relation $\entfde$
inherits this independence.

\begin{lemma}\label{lem:fourdecomp}
For all $A,B\in\For$,
\[
A\entfde B
\quad\Longleftrightarrow\quad
A\entfour{t}B\ \text{ and }\ A\entfour{1}B ;
\qquad\text{that is,}\qquad
\entfde\ =\ \entfour{t}\cap\entfour{1} .
\]
\end{lemma}

\begin{proof}
Fix a valuation $v$. In the two-element ontological base, the first
coordinate of $v(A)$ is at most that of $v(B)$ if and only if
$t\in v(A)\Rightarrow t\in v(B)$. In the epistemic base, the corresponding
condition is $1\in v(A)\Rightarrow 1\in v(B)$. Thus $v(A)\leqslant v(B)$
if and only if both implications hold. Quantifying over all valuations
gives the result by Definition~\ref{def:fdeent} and \eqref{eq:fourcomponents}.
\end{proof}

Lemma~\ref{lem:fourdecomp} gives a finite-valued instance of decomposition by
components. In the relational semantics of Definition~\ref{def:entailment},
$\entn=\entt\cap\entone$ holds by definition. The principal technical
contribution of the paper is the proof of the corresponding syntactic
decomposition, $\dern=\dert\cap\derone$, for the three separately axiomatised
calculi in Section~\ref{sec:algebra}.

The four-valued space validates the schemata
\begin{equation}\label{eq:paradox}
\negt A\wedge\nego A\entfde B,
\qquad
B\entfde\negt A\vee\nego A ,
\end{equation}
which are the onto-epistemic counterparts of \emph{ex falso quodlibet} and of
the law of excluded middle. Their validity is a consequence of the
bases being two-element: by
Lemma~\ref{lem:fourcomp}(1)--(3), for every value of $A$ the value of
$\negt A\wedge\nego A$ is $\langle f,0\rangle$ and that of
$\negt A\vee\nego A$ is $\langle t,1\rangle$.

The refinement adopted in \cite{bib51} consists in adjoining to each base a
third element representing indeterminacy. Since a value has two coordinates, indeterminacy
admits two forms, ontological and epistemic; they are written $u_t$ and $u_1$
accordingly. The bases become the three-element chains
\begin{equation}\label{eq:chains}
f\leqt u_t\leqt t,
\qquad
0\leqo u_1\leqo 1 ,
\end{equation}
and their product is the nine-element set
$\Nine=\{f,u_t,t\}\times\{0,u_1,1\}$.

\begin{definition}\label{def:nineorder}
For all $\langle s,u\rangle,\langle v,w\rangle\in\Nine$,
\[
\langle s,u\rangle\leqslant\langle v,w\rangle
\enskip\defiff\enskip
s\leqt v\ \text{ and }\ u\leqo w .
\]
\end{definition}

\begin{definition}\label{def:nineops}
The operations on $\Nine$ are given by \eqref{eq:fourops}, where now
$\sqcap_t,\sqcup_t$ and $\sqcap_1,\sqcup_1$ are the meet and join of the
respective chains in \eqref{eq:chains}, and $-_t$, $-_1$ interchange the two
extreme elements of the respective chain and fix its middle element.
\end{definition}

Definition~\ref{def:nineorder} makes $\Nine$ a lattice, denoted $\LatNine$,
displayed on the right of Figure~\ref{fig:lattices}. Being a product of two
finite chains, it is distributive, as is $\LatFour$.

\begin{definition}\label{def:nineent}
A valuation is a map $v\colon\Var\to\Nine$, extended to $\For$ as in
Definition~\ref{def:fourval}. For $A,B\in\For$,
\[
A\entu B
\enskip\defiff\enskip
v(A)\leqslant v(B)\ \text{ for every valuation } v .
\]
\end{definition}

The conditions in Lemma~\ref{lem:fourcomp} remain valid for $\Nine$. However, the
equivalences $t\in v(A)\iff f\notin v(A)$ and
$1\in v(A)\iff 0\notin v(A)$ no longer hold: a coordinate may instead take
the indeterminate value $u_t$ or $u_1$, respectively. By
Definition~\ref{def:nineops}, $u_t$ is a fixed point of $-_t$ and $u_1$ of
$-_1$. Consequently, if $v(A)=\langle u_t,u_1\rangle$, then $v(\negt
A)=v(\nego A)=v(A)$, so that $v(\negt A\wedge\nego A)=v(\negt A\vee\nego
A)=\langle u_t,u_1\rangle$; taking $v(B)$ below this value in the first case
and above it in the second refutes both schemata of \eqref{eq:paradox} with
respect to $\entu$. Indeterminacy of either kind thus turns a value into a
fixed point of the corresponding negation, and it is by this mechanism that the
schemata of \eqref{eq:paradox} cease to be valid.

Both $\LatFour$ and $\LatNine$ are distributive, being products of chains.
We now develop a relational semantics whose associated lattices need not be
distributive. The construction is based on Urquhart's representation theorem
for bounded lattices \cite{Urquhart} and its use by Allwein and
Dunn~\cite{AllweinDunn} in relational semantics.

In this representation, each lattice element is associated with a set of
maximal filter-ideal pairs whose filters contain it and with a corresponding
set whose ideals contain it. These sets are disjoint but need not cover the
carrier. A point outside both sets can be interpreted as a state in which the
represented element is neither established nor refuted.

The relational semantics below generalises the nine-valued semantics while
retaining its ontological and epistemic components. The values $u_t$ and $u_1$
are defined by nonmembership in both the truth and falsity sets of the
corresponding component. Such pointwise indeterminacy must be distinguished
from the existence of intermediate elements in the lattice of stable sets.
If a truth set and its corresponding falsity set do not cover the carrier,
the truth set must be nonempty and proper. Conversely, an intermediate stable
set need not give rise to pointwise indeterminacy. The relevant constructions
are introduced in Section~\ref{sec:dos}.

\subsection{Quadruply ordered sets}\label{sec:dos}

Urquhart's representation proceeds in two stages. The first is abstract: from a
set carrying two quasi-orders one extracts a bounded lattice. The second is
canonical: from a given nontrivial bounded lattice one constructs such a set,
and embeds the lattice into the bounded lattice extracted from it.

For the first stage, let $\langle X,\leqslant_1,\leqslant_2\rangle$ be a
\emph{doubly ordered set} in the sense used here, with $X$ non-empty
and $\leqslant_1,\leqslant_2$ quasi-orders. Urquhart's definition also
requires $x\leqslant_1 y$ and $x\leqslant_2 y$ to imply
$x=y$~\cite[Section~2]{Urquhart}. The stable-set construction works
for any pair of quasi-orders. Define
$\ell,r\colon\mathcal{P}(X)\to\mathcal{P}(X)$ by
\begin{equation}\label{eq:urqmaps}
\ell(Y)=\{x\mid\forall y\,(x\leqslant_1 y\Rightarrow y\notin Y)\},
\qquad
r(Y)=\{x\mid\forall y\,(x\leqslant_2 y\Rightarrow y\notin Y)\}.
\end{equation}
Restricted to increasing sets, $\ell$ and $r$ form an antitone Galois
connection; a set $Y$ is \emph{$\ell$-stable} if $\ell r(Y)=Y$ and
\emph{$r$-stable} if $r\ell(Y)=Y$. On the family $\Lstab(X)$ of $\ell$-stable
subsets of $X$ put
\begin{equation}\label{eq:urqops}
Y\wedge^{c}Z=Y\cap Z,
\qquad
Y\vee^{c}Z=\ell\bigl(r(Y)\cap r(Z)\bigr),
\qquad
0^{c}=\emptyset,
\qquad
1^{c}=X .
\end{equation}
Then $\langle\Lstab(X),\wedge^{c},\vee^{c},0^{c},1^{c}\rangle$ is a bounded
lattice~\cite{AllweinDunn}. The resulting lattice need not be distributive.
Intersections of $\ell$-stable sets are again $\ell$-stable, so meet is
set-theoretic, whereas unions need not be stable and join requires the closure
$\ell r$. These notions are stated in full, for the
duplicated case, in Definitions~\ref{def:galoismaps}--\ref{def:stable} below.

For the second stage, let $\Latt=\langle S,\wedge,\vee\rangle$ be a nontrivial
bounded lattice. A pair $(F,I)$ consisting of a filter $F$ and an ideal $I$ of $\Latt$
with $F\cap I=\emptyset$ is a \emph{filter-ideal pair}; it is \emph{maximal}
if neither component admits a proper extension disjoint from the other. Let
$X(\Latt)$ be the set of all maximal filter-ideal pairs, and for
$x\in X(\Latt)$ let $x_1,x_2$ denote its filter and ideal. Ordering
$X(\Latt)$ by
\begin{equation}\label{eq:urqorders}
x\leqslant_1 y\enskip\defiff\enskip x_1\subseteq y_1,
\qquad
x\leqslant_2 y\enskip\defiff\enskip x_2\subseteq y_2
\end{equation}
makes it a doubly ordered set in Urquhart's sense. If both inclusions
hold, maximality gives $x=y$. The map
\begin{equation}\label{eq:urqrep}
\beta\colon S\to\mathcal{P}(X(\Latt)),
\qquad
\beta a=\{x\in X(\Latt)\mid a\in x_1\}
\end{equation}
sends each element to an $\ell$-stable set and embeds $\Latt$ into
$\Lstab(X(\Latt))$; injectivity follows from the fact that a filter disjoint
from an ideal extends to a maximal pair separating them.

The corresponding set $r\beta a$ is described in terms of the ideal
components by
\begin{equation}\label{eq:urqsecond}
r\beta a=\{x\in X(\Latt)\mid a\in x_2\}.
\end{equation}
To verify \eqref{eq:urqsecond}, fix $x=(F,I)\in X(\Latt)$.
If $a\in I$, then every $(G,J)\in X(\Latt)$ with $I\subseteq J$ has
$a\in J$ and hence $a\notin G$. If $a\notin I$, the filter generated
by $a$ is disjoint from $I$. Extend this filter and $I$ to a maximal
disjoint pair $(G,J)$. Then $I\subseteq J$ and $a\in G$, which excludes
$x$ from $r\beta a$.

The sets $\beta a$ and $r\beta a$ are disjoint because $F\cap I=\emptyset$
for every maximal pair $(F,I)$. A point with $a$ in neither $F$ nor $I$
lies outside both sets. In the semantics below, a formula receives the
value $u_t$ or $u_1$ at points outside the truth and falsity sets of the
corresponding component.

In the Kripke semantics of Allwein and Dunn~\cite{AllweinDunn}, a doubly
ordered set is taken as primitive, and a valuation assigns to each
propositional variable a pair of subsets $(P_1,P_2)$ of $X$ subject to
\begin{equation}\label{eq:advaluation}
P_1=\ell(P_2),
\qquad
P_2=r(P_1).
\end{equation}
In this setting, \eqref{eq:advaluation} is
imposed as a condition on admissible valuations, rather than obtained from a
representation of a given lattice. Each component still determines the other.
Since the two sets are disjoint but need not cover $X$, they distinguish three
cases at a point: truth, falsity and indeterminacy. It is this three-valued
relational semantics that is extended below.

Write $x\forc{T}A$ for membership of $x$ in the first set associated with $A$
and $x\forc{F}A$ for membership in the second. The truth and falsity
conditions for conjunction and disjunction follow from the lattice operations in
\eqref{eq:urqops}. For disjunction, unfolding
$Y\vee^{c}Z=\ell(rY\cap rZ)$ gives a condition
quantified over $\leqslant_1$-successors,
\begin{equation}\label{eq:advee}
x\forc{T}A\vee B\iff\forall y\,\bigl(x\leqslant_1 y\Rightarrow
y\nforc{F}A\ \text{ or }\ y\nforc{F}B\bigr),
\end{equation}
whereas the corresponding falsity condition is pointwise, since
$r\bigl(\ell(rY\cap rZ)\bigr)=rY\cap rZ$:
\begin{equation}\label{eq:advee-f}
x\forc{F}A\vee B\iff x\forc{F}A\ \text{ and }\ x\forc{F}B .
\end{equation}
For conjunction, the truth condition is pointwise, since meet is intersection,
whereas the falsity condition quantifies over $\leqslant_2$-successors.

To interpret negation, we add an operator $\star\colon X\to X$ satisfying
\begin{equation}\label{eq:adstar}
x\leqslant_1 y\Rightarrow x^{\star}\leqslant_2 y^{\star},
\qquad
x\leqslant_2 y\Rightarrow x^{\star}\leqslant_1 y^{\star},
\qquad
x^{\star\star}=x.
\end{equation}
If $\Latt$ is equipped with a De Morgan negation $\mathord{\sim}$, the
canonical star is given by
$(x_1,x_2)^{\star}=(\mathord{\sim}x_2,\mathord{\sim}x_1)$, with
$\mathord{\sim}$ applied pointwise to subsets. The truth set of
$\mathord{\sim}A$ consists of the points $x$ such that
$x^{\star}\forc{F}A$. This gives
\begin{equation}\label{eq:adneg}
x\forc{T}\mathord{\sim}A\iff x^{\star}\forc{F}A,
\qquad
x\forc{F}\mathord{\sim}A\iff\forall y\,\bigl(x\leqslant_2
y\Rightarrow y^{\star}\nforc{F}A\bigr).
\end{equation}
The falsity condition is obtained by applying $r$ to this truth set. Conditions
\eqref{eq:adstar} ensure that the induced negation is an antitone involution
of the lattice $\Lstab(X)$ of $\ell$-stable sets.

The language $\Lang$ has two negations, each acting on its own component of
a value. We apply the construction above separately to the ontological and
epistemic components on a common carrier, using one pair of quasi-orders and
one involution for each component. From Allwein and Dunn we take the
doubly ordered structure, the maps $\ell$ and $r$, the condition
\eqref{eq:advaluation} and the star.

\begin{definition}\label{def:qos}
A \emph{quadruply ordered set} is a structure
$\Str=\langle X,\lt,\rt,\lo,\ro\rangle$ in which $X$ is a non-empty set, called
the \emph{carrier} of $\Str$, and each of $\lt,\rt,\lo,\ro$ is a quasi-order
on $X$.
\end{definition}

The four quasi-orders are grouped into the two pairs
$(\lt,\rt)$ and $(\lo,\ro)$, and each pair, taken together with $X$,
constitutes a doubly ordered set. Throughout, the
index $i\in\{t,1\}$ denotes the ontological or epistemic component,
respectively, and every notion introduced for a fixed $i$ applies to the
corresponding pair. No condition linking the
pairs is imposed at this stage; the sole such condition, which involves the two
star operations rather than the orders, appears in
Definition~\ref{def:frame} below.

Each pair of quasi-orders induces a pair of maps on $\mathcal{P}(X)$.

\begin{definition}\label{def:galoismaps}
Let $\Str$ be a quadruply ordered set and $Y\subseteq X$. Put
\begin{align}
\lmap{i}(Y)&=\{x\in X\mid \forall y\in X\,(x\ordl{i}y\Rightarrow y\notin Y)\},
  \label{eq:lmap}\\
\rmap{i}(Y)&=\{x\in X\mid \forall y\in X\,(x\ordr{i}y\Rightarrow y\notin Y)\}.
  \label{eq:rmap}
\end{align}
\end{definition}

\begin{definition}\label{def:increasing}
A set $Y\subseteq X$ is \emph{$\ordl{i}$-increasing} if $x\in Y$ and
$x\ordl{i}y$ jointly imply $y\in Y$, for all $x,y\in X$; it is
\emph{$\ordr{i}$-increasing} if $x\in Y$ and $x\ordr{i}y$ jointly imply $y\in
Y$.
\end{definition}

\begin{definition}\label{def:stable}
A set $Y\subseteq X$ is \emph{$\lmap{i}$-stable} if $\lmap{i}\rmap{i}(Y)=Y$,
and \emph{$\rmap{i}$-stable} if $\rmap{i}\lmap{i}(Y)=Y$.
\end{definition}

The maps $\lmap{i}$ and $\rmap{i}$ are antitone
(Lemma~\ref{lem:increasing}). On $\mathcal{P}(X)$, the Galois equivalence
can fail. For example, take $X=\{x,y\}$ with $x\ordl{i}y$ the only
nonreflexive comparison and let $\ordr{i}$ be equality.
For $Y=\{x\}$ and $Z=\{y\}$, we have $Z=\rmap{i}(Y)$ but
$\lmap{i}(Z)=\emptyset$, so $Y\nsubseteq\lmap{i}(Z)$. Here $Y$
fails to be $\ordl{i}$-increasing. Lemma~\ref{lem:galois} gives the
connection on the corresponding families of increasing sets.

\begin{lemma}\label{lem:increasing}
Let $i\in\{t,1\}$. For every $Y\subseteq X$ the set $\lmap{i}(Y)$ is
$\ordl{i}$-increasing and the set $\rmap{i}(Y)$ is $\ordr{i}$-increasing;
moreover both maps are antitone, that is, $Y\subseteq Z$ implies
$\lmap{i}(Z)\subseteq\lmap{i}(Y)$ and $\rmap{i}(Z)\subseteq\rmap{i}(Y)$.
\end{lemma}

\begin{proof}
Let $x\in\lmap{i}(Y)$ and $x\ordl{i}x'$. If $x'\ordl{i}y$, then $x\ordl{i}y$ by
transitivity, whence $y\notin Y$. Thus $x'\in\lmap{i}(Y)$. The argument for
$\rmap{i}(Y)$ is the same. Antitonicity follows directly from the definitions.
\end{proof}

\begin{lemma}\label{lem:galois}
Let $i\in\{t,1\}$ and let $Y,Z\subseteq X$ be $\ordl{i}$-increasing and
$\ordr{i}$-increasing, respectively. Then
\begin{equation}\label{eq:galois}
Y\subseteq\lmap{i}(Z)
\quad\Longleftrightarrow\quad
Y\cap Z=\emptyset
\quad\Longleftrightarrow\quad
Z\subseteq\rmap{i}(Y).
\end{equation}
Hence $\lmap{i}$ and $\rmap{i}$, restricted to $\ordr{i}$-increasing and
$\ordl{i}$-increasing subsets of $X$, respectively, form an antitone Galois
connection between these two families ordered by inclusion.
\end{lemma}

\begin{proof}
Suppose $Y\cap Z=\emptyset$. If $x\in Y$ and $x\ordl{i}z$, then $z\in Y$
because $Y$ is $\ordl{i}$-increasing, and hence $z\notin Z$. Thus
$Y\subseteq\lmap{i}(Z)$. Conversely, this inclusion implies
$Y\cap Z=\emptyset$ by reflexivity of $\ordl{i}$. The second equivalence
follows in the same way, using that $Z$ is $\ordr{i}$-increasing.
\end{proof}

\begin{corollary}\label{cor:closure}
Let $i\in\{t,1\}$. Then
\begin{enumerate}
\item[(1)] $\lmap{i}\rmap{i}$ is a closure operator on the family of
$\ordl{i}$-increasing subsets of $X$, and $\rmap{i}\lmap{i}$ is a closure
operator on the family of $\ordr{i}$-increasing subsets of $X$;
\item[(2)] the $\lmap{i}$-stable sets are closed under arbitrary
intersections, and so are the $\rmap{i}$-stable sets.
\end{enumerate}
\end{corollary}

\begin{proof}
Item (1) follows from the antitone Galois connection of
Lemma~\ref{lem:galois}. By Lemma~\ref{lem:increasing}, every
$\lmap{i}$-stable set is $\ordl{i}$-increasing, and every $\rmap{i}$-stable
set is $\ordr{i}$-increasing. Thus the stable sets are precisely the fixed
points of the corresponding closure operators. Both families of increasing
sets are closed under arbitrary intersections, so the same is true of the
fixed points of their closure operators, giving (2).
\end{proof}

In each component, the truth and falsity sets of a formula will form a pair
satisfying the following condition.

\begin{definition}\label{def:conjugate}
A pair $(Y,Z)$ of subsets of $X$ is \emph{$i$-conjugate} if
\begin{equation}\label{eq:conjugate}
\rmap{i}(Y)=Z
\qquad\text{and}\qquad
\lmap{i}(Z)=Y .
\end{equation}
\end{definition}

\begin{lemma}\label{lem:conjugate}
Let $(Y,Z)$ be $i$-conjugate. Then $Y$ is $\lmap{i}$-stable and
$\ordl{i}$-increasing, $Z$ is $\rmap{i}$-stable and $\ordr{i}$-increasing,
and $Y\cap Z=\emptyset$. Conversely, if $Y$ is $\lmap{i}$-stable, then
$(Y,\rmap{i}(Y))$ is $i$-conjugate, and if $Z$ is $\rmap{i}$-stable, then
$(\lmap{i}(Z),Z)$ is $i$-conjugate. Thus $\rmap{i}$ and $\lmap{i}$ restrict
to mutually inverse bijections between the $\lmap{i}$-stable and
$\rmap{i}$-stable subsets of $X$.
\end{lemma}

\begin{proof}
Substituting the first equation of \eqref{eq:conjugate} into the second gives
$\lmap{i}\rmap{i}(Y)=Y$; substituting the second into the first gives
$\rmap{i}\lmap{i}(Z)=Z$. As $Y$ and $Z$ are values of $\lmap{i}$ and $\rmap{i}$
respectively, Lemma~\ref{lem:increasing} shows that $Y$ is
$\ordl{i}$-increasing and $Z$ is $\ordr{i}$-increasing.
Lemma~\ref{lem:galois} then yields $Y\cap Z=\emptyset$.
For the converse, let $\lmap{i}\rmap{i}(Y)=Y$ and put $Z:=\rmap{i}(Y)$: the
first equation of \eqref{eq:conjugate} holds by the choice of $Z$, the second
by the assumption. The remaining case is symmetric. The equations
\eqref{eq:conjugate} also show that the two restrictions are mutually inverse.
\end{proof}

\subsection[NDst\_oe-frames]{$NDst_{oe}$-frames}\label{sec:frames}

To interpret the two negations, we equip a quadruply ordered set with two
commuting involutions.

\begin{definition}\label{def:frame}
An \emph{$NDst_{oe}$-frame} is a structure
\[
\Frm=\langle X,\lt,\rt,\lo,\ro,\Nt,\No\rangle
\]
in which
$\langle X,\lt,\rt,\lo,\ro\rangle$ is a quadruply ordered set and $\Nt,\No$
are unary operations on $X$ satisfying, for all $x,y\in X$ and each
$i\in\{t,1\}$,
\begin{align}
N_{i}N_{i}x&=x, \label{eq:frinv}\\
x\ordl{i}y&\iff N_{i}x\ordr{i}N_{i}y, \label{eq:frswap}\\
\Nt\No x&=\No\Nt x. \label{eq:frcomm}
\end{align}
\end{definition}

Conditions \eqref{eq:frinv} and \eqref{eq:frswap} apply to each component
separately; condition \eqref{eq:frcomm} relates the two involutions.

\begin{lemma}\label{lem:frame}
Let $\Frm$ be an $NDst_{oe}$-frame, $x,y\in X$ and $i\in\{t,1\}$. Then
\begin{align}
x\ordr{i}y&\iff N_{i}x\ordl{i}N_{i}y, \label{eq:frswapback}\\
N_{i}x\ordl{i}y&\iff x\ordr{i}N_{i}y. \label{eq:fradj}
\end{align}
\end{lemma}

\begin{proof}
For \eqref{eq:frswapback}, instantiate \eqref{eq:frswap} at $N_{i}x$ and
$N_{i}y$: this gives
$N_{i}x\ordl{i}N_{i}y\iff N_{i}N_{i}x\ordr{i}N_{i}N_{i}y$, and the right-hand
side is $x\ordr{i}y$ by \eqref{eq:frinv}. For \eqref{eq:fradj}, instantiate
\eqref{eq:frswap} at $N_{i}x$ and $y$: this gives
$N_{i}x\ordl{i}y\iff N_{i}N_{i}x\ordr{i}N_{i}y$, and the left argument on the
right-hand side is $x$ by \eqref{eq:frinv}.
\end{proof}

Each $N_{i}$ is an isomorphism from $\langle X,\ordl{i}\rangle$ onto
$\langle X,\ordr{i}\rangle$, with inverse $N_{i}$. For each $i$, conditions
\eqref{eq:frinv} and \eqref{eq:frswap} are equivalent to the star conditions
\eqref{eq:adstar}.

\begin{remark}\label{rem:commindep}
Condition \eqref{eq:frcomm} is independent of \eqref{eq:frinv} and
\eqref{eq:frswap}. Let $X=\{a,b,c\}$, let all four quasi-orders be equality,
and define
\[
\begin{aligned}
\Nt a&=b, & \Nt b&=a, & \Nt c&=c,\\
\No a&=a, & \No b&=c, & \No c&=b.
\end{aligned}
\]
Then \eqref{eq:frinv} and \eqref{eq:frswap} hold, but
\[
\Nt\No a=b\ne c=\No\Nt a.
\]
\end{remark}

\subsection{Valuations and truth conditions}\label{sec:valuations}
To interpret the formulae of $\Lang$ in a frame we first fix how values are
assigned to propositional variables. Each variable receives one pair of sets in
each component, the first member indicating where it is true and the second
where it is false, and each pair is required to be conjugate in the sense of
Definition~\ref{def:conjugate}.

\begin{definition}\label{def:valuation}
Let $\Frm$ be an $NDst_{oe}$-frame with carrier $X$. A \emph{valuation} $V$
on $\Frm$ assigns to each $p\in\Var$ two pairs of subsets of $X$: a
$t$-conjugate pair $(\Pset{p}{t},\Pset{p}{f})$ and a $1$-conjugate pair
$(\Pset{p}{1},\Pset{p}{0})$. For $p\in\Var$, $x\in X$ and
$j\in\{t,f,1,0\}$ put
\[
x\forc{j}p\enskip\defiff\enskip x\in\Pset{p}{j}.
\]
A \emph{model} is a frame together with a valuation on it.
\end{definition}

By Lemma~\ref{lem:conjugate} the members of a conjugate pair are disjoint.
Their union, however, need not be $X$. At a point outside both sets, $p$ is
neither true nor false in the corresponding component. The indeterminate
values are defined in Definition~\ref{def:indeterminacy}.

We now extend these relations to all formulae.

\begin{definition}\label{def:forc}
In a model with carrier $X$, let $x\in X$ and
$A,B\in\For$. The relations $\forc{t},\forc{f},\forc{1},\forc{0}$ are
defined on variables by Definition~\ref{def:valuation} and extended
recursively to compound formulae. The pointwise conditions for the
binary connectives are
\begin{align}
x\forc{t}A\wedge B&\iff x\forc{t}A\ \text{ and }\ x\forc{t}B,\label{eq:cl1}\\
x\forc{1}A\wedge B&\iff x\forc{1}A\ \text{ and }\ x\forc{1}B,\label{eq:cl2}\\
x\forc{f}A\vee B&\iff x\forc{f}A\ \text{ and }\ x\forc{f}B,\label{eq:cl3}\\
x\forc{0}A\vee B&\iff x\forc{0}A\ \text{ and }\ x\forc{0}B.\label{eq:cl4}
\end{align}
The corresponding quantified conditions are
\begin{align}
x\forc{t}A\vee B&\iff\forall y\,\bigl(x\lt y\Rightarrow
  y\nforc{f}A\ \text{ or }\ y\nforc{f}B\bigr),\label{eq:cl5}\\
x\forc{1}A\vee B&\iff\forall y\,\bigl(x\lo y\Rightarrow
  y\nforc{0}A\ \text{ or }\ y\nforc{0}B\bigr),\label{eq:cl6}\\
x\forc{f}A\wedge B&\iff\forall y\,\bigl(x\rt y\Rightarrow
  y\nforc{t}A\ \text{ or }\ y\nforc{t}B\bigr),\label{eq:cl7}\\
x\forc{0}A\wedge B&\iff\forall y\,\bigl(x\ro y\Rightarrow
  y\nforc{1}A\ \text{ or }\ y\nforc{1}B\bigr).\label{eq:cl8}
\end{align}
For each negation in its own component, put
\begin{align}
x\forc{t}\negt A&\iff \Nt x\forc{f}A,\label{eq:cl9}\\
x\forc{f}\negt A&\iff\forall y\,\bigl(x\rt y\Rightarrow
  \Nt y\nforc{f}A\bigr),\label{eq:cl10}\\
x\forc{1}\nego A&\iff \No x\forc{0}A,\label{eq:cl11}\\
x\forc{0}\nego A&\iff\forall y\,\bigl(x\ro y\Rightarrow
  \No y\nforc{0}A\bigr).\label{eq:cl12}
\end{align}
The \emph{transparency} conditions say that each negation leaves the
other component unchanged.
\begin{align}
x\forc{t}\nego A&\iff x\forc{t}A,
&x\forc{f}\nego A&\iff x\forc{f}A,\label{eq:cl13}\\
x\forc{1}\negt A&\iff x\forc{1}A,
&x\forc{0}\negt A&\iff x\forc{0}A.\label{eq:cl14}
\end{align}
\end{definition}

For $A\in\For$ and $j\in\{t,f,1,0\}$ write
$\Pset{A}{j}=\{x\in X\mid x\forc{j}A\}$, extending the notation of
Definition~\ref{def:valuation}.

\begin{definition}\label{def:indeterminacy}
For every $A\in\For$ and every $x\in X$ put
\[
x\forc{u_{t}}A\enskip\defiff\enskip x\nforc{t}A\ \text{ and }\ x\nforc{f}A,
\qquad
x\forc{u_{1}}A\enskip\defiff\enskip x\nforc{1}A\ \text{ and }\ x\nforc{0}A .
\]
\end{definition}

For $Y\subseteq X$ and $i\in\{t,1\}$ write $N_{i}(Y):=\{N_{i}x\mid x\in Y\}$,
which by \eqref{eq:frinv} coincides with $\{x\mid N_{i}x\in Y\}$.

\begin{lemma}\label{lem:starmaps}
For every $Y\subseteq X$ and every $i\in\{t,1\}$,
\begin{equation}\label{eq:starmaps}
\rmap{i}\bigl(N_{i}(Y)\bigr)=N_{i}\bigl(\lmap{i}(Y)\bigr),
\qquad
\lmap{i}\bigl(N_{i}(Y)\bigr)=N_{i}\bigl(\rmap{i}(Y)\bigr).
\end{equation}
\end{lemma}

\begin{proof}
By \eqref{eq:rmap} and \eqref{eq:lmap}, for every $x\in X$,
\[
\begin{aligned}
x\in\rmap{i}(N_{i}(Y))
&\iff\forall y\,(x\ordr{i}y\Rightarrow N_{i}y\notin Y)\\
&\iff\forall z\,(N_{i}x\ordl{i}z\Rightarrow z\notin Y)\\
&\iff N_{i}x\in\lmap{i}(Y)\\
&\iff x\in N_{i}(\lmap{i}(Y)).
\end{aligned}
\]
The second equivalence uses the substitution $z=N_{i}y$, which is reversible
by \eqref{eq:frinv}, together with \eqref{eq:frswapback}. This proves the
first identity. Applying it to $N_{i}(Y)$ and using \eqref{eq:frinv} gives
\[
\rmap{i}(Y)=N_{i}\bigl(\lmap{i}(N_{i}(Y))\bigr).
\]
Applying $N_{i}$ to both sides and using involutivity once more yields the
second identity.
\end{proof}

The following lemma shows that the conjugacy conditions of
Definition~\ref{def:valuation} are preserved under the semantic conditions in
Definition~\ref{def:forc}.

\begin{lemma}\label{lem:conjugacypreserved}
In every $NDst_{oe}$-model and for every $A\in\For$, the pair
$(\Pset{A}{t},\Pset{A}{f})$ is $t$-conjugate and $(\Pset{A}{1},\Pset{A}{0})$ is
$1$-conjugate; that is,
\begin{equation}\label{eq:conjpres}
\rmap{t}(\Pset{A}{t})=\Pset{A}{f},\quad
\lmap{t}(\Pset{A}{f})=\Pset{A}{t},\quad
\rmap{1}(\Pset{A}{1})=\Pset{A}{0},\quad
\lmap{1}(\Pset{A}{0})=\Pset{A}{1}.
\end{equation}
\end{lemma}

\begin{proof}
Fix a model and $i\in\{t,1\}$. Write $T_A$ and $F_A$ for the truth and
falsity sets of $A$ in component $i$. We show by structural induction
that $(T_A,F_A)$ is $i$-conjugate. For variables this is
Definition~\ref{def:valuation}.

The semantic conditions for conjunction give
$T_{A\wedge B}=T_A\cap T_B$ and $F_{A\wedge B}=\rmap{i}(T_{A\wedge B})$.
The induction hypothesis and Corollary~\ref{cor:closure}(2) give
$\lmap{i}$-stability of this intersection, so the pair is conjugate by
Lemma~\ref{lem:conjugate}. For disjunction, $F_{A\vee B}=F_A\cap F_B$
and $T_{A\vee B}=\lmap{i}(F_{A\vee B})$. Here conjugacy follows from the
$\rmap{i}$-stability of the intersection.

For the negation $\neg_i$ of component $i$, the semantic conditions
give $T_{\neg_i A}=N_i(F_A)$ and
$F_{\neg_i A}=\rmap{i}(T_{\neg_i A})$. Lemma~\ref{lem:starmaps} and
the induction hypothesis yield
\[
F_{\neg_i A}=N_i(\lmap{i}(F_A))=N_i(T_A),
\qquad
\lmap{i}(F_{\neg_i A})=N_i(\rmap{i}(T_A))=T_{\neg_i A}.
\]
The other negation leaves both sets unchanged in component $i$, so
the induction hypothesis applies.
\end{proof}

\begin{corollary}\label{cor:stableformulae}
For every $A\in\For$, the set $\Pset{A}{t}$ is $\lmap{t}$-stable,
$\Pset{A}{f}$ is $\rmap{t}$-stable, $\Pset{A}{1}$ is $\lmap{1}$-stable and
$\Pset{A}{0}$ is $\rmap{1}$-stable; moreover $\Pset{A}{t}\cap\Pset{A}{f}
=\emptyset$ and $\Pset{A}{1}\cap\Pset{A}{0}=\emptyset$.
\end{corollary}

\begin{proof}
Immediate from Lemma~\ref{lem:conjugacypreserved} and
Lemma~\ref{lem:conjugate}.
\end{proof}

\begin{corollary}\label{cor:threevalued}
For every $A\in\For$ and every $x\in X$, exactly one of $x\forc{t}A$,
$x\forc{u_{t}}A$, $x\forc{f}A$ holds, and exactly one of $x\forc{1}A$,
$x\forc{u_{1}}A$, $x\forc{0}A$. The ontological and epistemic components
therefore determine total maps from $\For\times X$ to $\{t,u_{t},f\}$ and
$\{1,u_{1},0\}$, respectively.
\end{corollary}

\begin{proof}
At most one of the three holds: $\Pset{A}{t}\cap\Pset{A}{f}=\emptyset$ by
Corollary~\ref{cor:stableformulae}, and $x\forc{u_{t}}A$ excludes the other two
by Definition~\ref{def:indeterminacy}. At least one holds: if neither
$x\forc{t}A$ nor $x\forc{f}A$, then $x\forc{u_{t}}A$ by the same definition.
The epistemic case is identical.
\end{proof}

For $i\in\{t,1\}$ write $\Lstab_{i}$ for the family of
$\lmap{i}$-stable subsets of $X$. The closure operator of
Corollary~\ref{cor:closure} gives a bounded lattice under inclusion with
the operations in \eqref{eq:urqops} for the corresponding pair of
quasi-orders. In particular, the join is $\lmap{i}\rmap{i}(Y\cup Z)$,
and $\rmap{i}(Y\cup Z)=\rmap{i}(Y)\cap\rmap{i}(Z)$ by definition.

\begin{lemma}\label{lem:starinvolution}
Let $\Frm$ be an $NDst_{oe}$-frame with carrier $X$ and let
$i\in\{t,1\}$. The map
\[
\Nc{i}(Y):=N_{i}(\rmap{i}(Y))
\]
is an antitone involution of the lattice $\Lstab_{i}$.
\end{lemma}

\begin{proof}
Let $Y$ be $\lmap{i}$-stable, so that $\lmap{i}\rmap{i}Y=Y$. By
Lemma~\ref{lem:starmaps},
\[
\lmap{i}\rmap{i}\bigl(N_{i}(\rmap{i}Y)\bigr)
=\lmap{i}\bigl(N_{i}(\lmap{i}\rmap{i}Y)\bigr)
=\lmap{i}(N_{i}(Y))
=N_{i}(\rmap{i}Y).
\]
Thus $\Nc{i}$ maps $\Lstab_{i}$ into itself. It is antitone because
$\rmap{i}$ is antitone by Lemma~\ref{lem:increasing} and
$N_{i}(\,\cdot\,)$ preserves inclusion. Applying Lemma~\ref{lem:starmaps}
once more and using \eqref{eq:frinv},
\[
\Nc{i}\Nc{i}(Y)
=N_{i}\bigl(\rmap{i}(N_{i}(\rmap{i}Y))\bigr)
=N_{i}\bigl(N_{i}(\lmap{i}\rmap{i}Y)\bigr)
=\lmap{i}\rmap{i}Y=Y.
\]
\end{proof}

\begin{lemma}\label{lem:negcomponents}
Let a model on an $NDst_{oe}$-frame with carrier $X$ be given and let
$A\in\For$. Then
\begin{enumerate}
\item[(1)] $\Pset{\nego A}{t}=\Pset{A}{t}$, $\Pset{\nego A}{f}=\Pset{A}{f}$,
$\Pset{\negt A}{1}=\Pset{A}{1}$ and $\Pset{\negt A}{0}=\Pset{A}{0}$;
\item[(2)] $\Pset{\negt A}{t}=N_{t}(\Pset{A}{f})$ and
$\Pset{\negt A}{f}=N_{t}(\Pset{A}{t})$; likewise
$\Pset{\nego A}{1}=N_{1}(\Pset{A}{0})$ and
$\Pset{\nego A}{0}=N_{1}(\Pset{A}{1})$;
\item[(3)] $\Nc{t}(\Pset{A}{t})=\Pset{\negt A}{t}$ and
$\Nc{1}(\Pset{A}{1})=\Pset{\nego A}{1}$.
\end{enumerate}
\end{lemma}

\begin{proof}
Item (1) follows from the transparency conditions \eqref{eq:cl13} and
\eqref{eq:cl14}.

The induction step for $\neg_i$ in Lemma~\ref{lem:conjugacypreserved}
gives (2) for $i=t$ and $i=1$.

For (3), $\Nc{t}(\Pset{A}{t})=N_{t}(\rmap{t}\Pset{A}{t})
=N_{t}(\Pset{A}{f})=\Pset{\negt A}{t}$ by
Lemma~\ref{lem:conjugacypreserved} and item (2). The epistemic identity is
obtained in the same way.
\end{proof}

The semantic conditions of Definition~\ref{def:forc} can now be expressed
using the lattice operations and the antitone involutions $\Nc{t},\Nc{1}$
defined in Lemma~\ref{lem:starinvolution}.

\begin{corollary}\label{cor:homomorphism}
For all $A,B\in\For$,
\[
\Pset{A\wedge B}{t}=\Pset{A}{t}\cap\Pset{B}{t},
\qquad
\Pset{A\vee B}{t}=\lmap{t}\bigl(\rmap{t}\Pset{A}{t}\cap\rmap{t}\Pset{B}{t}\bigr),
\qquad
\Pset{\negt A}{t}=\Nc{t}\bigl(\Pset{A}{t}\bigr),
\qquad
\Pset{\nego A}{t}=\Pset{A}{t}.
\]
The value $\Pset{A}{t}$ is thus computed from the values of the variables
occurring in $A$ by the operations \eqref{eq:urqops} of $\Lstab_{t}$, with
$\negt$ read as $\Nc{t}$ and $\nego$ as the identity. In particular
\begin{equation}\label{eq:demorgan}
\Pset{\negt(A\wedge B)}{t}=\Pset{\negt A\vee\negt B}{t},
\qquad
\Pset{\negt(A\vee B)}{t}=\Pset{\negt A\wedge\negt B}{t}.
\end{equation}
Mirror statements hold for $\Pset{A}{1}$ in $\Lstab_{1}$, with $\nego$ read as
$\Nc{1}$ and $\negt$ as the identity.
\end{corollary}

\begin{proof}
The first identity is \eqref{eq:cl1}; the second is \eqref{eq:cl5} together
with Lemma~\ref{lem:conjugacypreserved}, by which $\Pset{A}{f}=
\rmap{t}\Pset{A}{t}$; the third and the fourth are
Lemma~\ref{lem:negcomponents}(3) and~(1). Since $\Nc{t}$ is an antitone
involution of $\Lstab_{t}$ by Lemma~\ref{lem:starinvolution}, it interchanges
meets and joins, yielding \eqref{eq:demorgan}.
\end{proof}

\subsection{Three entailment relations}\label{sec:entailment}

A two-component semantics admits more than one notion of entailment, since
truth may be required to be preserved in one component, in the other, or in
both.

\begin{definition}\label{def:entailment}
Let $A,B\in\For$. The relations $A\entt B$ and $A\entone B$ hold if and
only if the respective conditions
\begin{align}
x\forc{t}A&\Rightarrow x\forc{t}B,\label{eq:entt}\\
x\forc{1}A&\Rightarrow x\forc{1}B\label{eq:entone}
\end{align}
hold at every point $x$ of every model. Joint entailment is defined by
\begin{equation}\label{eq:entn}
A\entn B\enskip\defiff\enskip A\entt B\ \text{ and }\ A\entone B .
\end{equation}
\end{definition}

The following lemma shows that, over all points of a fixed model, forward
preservation of truth is equivalent to backward preservation of falsity in
each component.

\begin{lemma}\label{lem:reduction}
Let $A,B\in\For$ and let a model on an $NDst_{oe}$-frame with carrier $X$ be
given. Then
\begin{equation}\label{eq:reduction}
\Pset{A}{t}\subseteq\Pset{B}{t}
\iff
\Pset{B}{f}\subseteq\Pset{A}{f},
\qquad
\Pset{A}{1}\subseteq\Pset{B}{1}
\iff
\Pset{B}{0}\subseteq\Pset{A}{0}.
\end{equation}
\end{lemma}

\begin{proof}
The maps $\rmap{t}$ and $\lmap{t}$ are antitone by
Lemma~\ref{lem:increasing}. Assume $\Pset{A}{t}\subseteq\Pset{B}{t}$.
Applying $\rmap{t}$ gives
$\rmap{t}(\Pset{B}{t})\subseteq\rmap{t}(\Pset{A}{t})$, and by
Lemma~\ref{lem:conjugacypreserved} the two sides are $\Pset{B}{f}$ and
$\Pset{A}{f}$. Conversely, assume $\Pset{B}{f}\subseteq\Pset{A}{f}$. Applying
$\lmap{t}$ gives $\lmap{t}(\Pset{A}{f})\subseteq\lmap{t}(\Pset{B}{f})$, and by
Lemma~\ref{lem:conjugacypreserved} again the two sides are $\Pset{A}{t}$ and
$\Pset{B}{t}$. The epistemic case replaces $\rmap{t},\lmap{t}$ by
$\rmap{1},\lmap{1}$ throughout.
\end{proof}

By Definition~\ref{def:entailment}, the third relation is the intersection
of the two component relations:
\begin{equation}\label{eq:semdecomp}
\entn\ =\ \entt\cap\entone .
\end{equation}
The corresponding equality for the independently axiomatised derivability
relations requires a separate proof.

\subsection{The finite-valued semantics as special cases}\label{sec:onepoint}

We recover the four- and nine-valued semantics of Section~\ref{sec:finite}
on a one-point frame and a three-point frame, respectively. Their stable-set
lattices reproduce the component bases, and the finite value of every formula
is realised at a distinguished point $c$.

\begin{samepage}
\begin{definition}\label{def:finiteframes}
Define two structures $\Frm_{2}$ and $\Frm_{3}$ as follows. The subscripts
indicate the number of values in each component.
\begin{enumerate}
\item[(1)] $\Frm_{2}$ has carrier $X_{2}=\{c\}$, all four quasi-orders are
the equality relation, and both involutions $\Nt,\No$ are the identity.
\item[(2)] $\Frm_{3}$ has carrier $X_{3}=\{a,b,c\}$. For each $i\in\{t,1\}$,
the quasi-orders $\ordl{i}$ and $\ordr{i}$ are the linear orders determined by
\[
c\ordl{i}b\ordl{i}a,
\qquad
c\ordr{i}a\ordr{i}b,
\]
and the involutions satisfy
\[
N_{i}a=b,\qquad N_{i}b=a,\qquad N_{i}c=c.
\]
\end{enumerate}
\end{definition}
\end{samepage}

\begin{lemma}\label{lem:finiteframes}
Both structures are $NDst_{oe}$-frames. For each $i\in\{t,1\}$ the
families of $\lmap{i}$-stable sets are
\[
\{\emptyset,X_{2}\}\quad\text{on }\Frm_{2},
\qquad
\{\emptyset,\{a\},X_{3}\}\quad\text{on }\Frm_{3}.
\]
In both lattices, meet is intersection and join is union. The antitone
involution $\Nc{i}$ interchanges the empty set and the carrier and, on
$\Frm_{3}$, fixes $\{a\}$.
\end{lemma}

\begin{proof}
The frame conditions for $\Frm_{2}$ are immediate. On its carrier,
$\lmap{i}$ and $\rmap{i}$ interchange $\emptyset$ and $X_{2}$, giving the
stated stable sets and involution.

On $\Frm_{3}$, the map exchanging $a$ and $b$ and fixing $c$ is an involutive
isomorphism from $\langle X_{3},\ordl{i}\rangle$ onto
$\langle X_{3},\ordr{i}\rangle$. Thus \eqref{eq:frinv} and
\eqref{eq:frswap} hold; \eqref{eq:frcomm} holds because the two involutions
coincide. By Lemma~\ref{lem:increasing}, every $\lmap{i}$-stable set is
$\ordl{i}$-increasing, so it is one of
$\emptyset,\{a\},\{a,b\},X_{3}$. Directly from \eqref{eq:lmap} and
\eqref{eq:rmap},
\[
\rmap{i}(\{a\})=\{b\},\qquad
\lmap{i}(\{b\})=\{a\},\qquad
\rmap{i}(\{a,b\})=\emptyset.
\]
Since $\lmap{i}(\emptyset)=\rmap{i}(\emptyset)=X_{3}$ and
$\lmap{i}(X_{3})=\rmap{i}(X_{3})=\emptyset$, the empty set and the carrier
are stable. Also,
$\lmap{i}\rmap{i}(\{a,b\})=X_{3}\ne\{a,b\}$, so $\{a,b\}$ is not
$\lmap{i}$-stable. Thus the stable sets are exactly
$\emptyset,\{a\},X_{3}$. Moreover
$\Nc{i}(\{a\})=N_{i}(\{b\})=\{a\}$, and $\Nc{i}$ interchanges the two
extreme sets. Each displayed family is a chain under inclusion, so meet and
join are intersection and union.
\end{proof}

Recall that the truth and falsity sets assigned to variables form conjugate
pairs by Definition~\ref{def:valuation}, and
Lemma~\ref{lem:conjugacypreserved} establishes the same property for all
formulae. On $\Frm_{3}$, the possible ontological pairs
$(\Pset{A}{t},\Pset{A}{f})$ are listed below, together with the finite
component values they represent:
\[
\begin{array}{c|cc}
\text{value}&\text{truth set}&\text{falsity set}\\
\hline
f&\emptyset&X_{3}\\
u_t&\{a\}&\{b\}\\
t&X_{3}&\emptyset
\end{array}
\]
For the epistemic component, replace $f,u_t,t$ by $0,u_1,1$, respectively.
On $\Frm_{2}$, use the first and last rows with $X_{2}$ in place of $X_{3}$.
The two components are chosen independently.
At $c$, the first and last rows give falsity and truth, respectively. On
$\Frm_{3}$, the middle row places $c$ outside both sets, giving indeterminacy
by Definition~\ref{def:indeterminacy}.

\begin{lemma}\label{lem:finitevaluations}
The correspondence above induces a bijection between valuations on
$\Frm_{2}$ and $\Four$-valued valuations, and between valuations on
$\Frm_{3}$ and $\Nine$-valued valuations. For corresponding relational
valuations $V$ (Definition~\ref{def:valuation}) and finite-valued valuations
$v$, evaluating any formula $A$ in the relational semantics yields exactly
the two conjugate pairs representing $v(A)$ in the table.
\end{lemma}

\begin{proof}
By Lemma~\ref{lem:conjugate}, an admissible assignment to a variable is
determined by its two stable truth sets. By Lemma~\ref{lem:finiteframes},
in each component the table gives an order isomorphism from the chain of
values to the lattice of stable truth sets, using the extreme rows on
$\Frm_2$. This gives the stated bijections of valuations.

We extend the correspondence to formulae by structural induction. The
variable case holds by construction. For the inductive steps, we check
that the table preserves each operation. It suffices to track truth sets,
since conjugacy uniquely determines the corresponding falsity sets. On
$\Frm_{3}$, let $h$ assign to each ontological component value its truth
set in the table:
\[
h(f)=\emptyset,\qquad h(u_t)=\{a\},\qquad h(t)=X_{3}.
\]
\begin{samepage}
Suppose the correspondence holds for $B$ and $C$, and let $\alpha,\beta$ be
the ontological coordinates of $v(B),v(C)$, respectively. Since $h$
preserves order, it takes minimum to intersection. By
Corollary~\ref{cor:homomorphism},
\[
\Pset{B\wedge C}{t}
=\Pset{B}{t}\cap\Pset{C}{t}
=h(\alpha)\cap h(\beta)
=h\bigl(\min(\alpha,\beta)\bigr).
\]
The last set represents the ontological coordinate of $v(B\wedge C)$.\par
\end{samepage}

\begin{samepage}
For disjunction the same argument uses maximum and union, by
Lemma~\ref{lem:finiteframes}. For $\negt$, Corollary~\ref{cor:homomorphism}
uses $\Nc{t}$, which interchanges the extreme sets and fixes the middle
set by Lemma~\ref{lem:finiteframes}; $\nego$ leaves the truth set
unchanged. These agree with the finite operations
of \eqref{eq:fourops} and Definition~\ref{def:nineops}. The epistemic
component is analogous, with the roles of the two negations interchanged.
On $\Frm_{2}$ the argument uses only the extreme rows of the table. This
completes the induction for truth sets. By
Lemma~\ref{lem:conjugacypreserved}, their falsity sets are the conjugate
sets, hence are also given by the table.\qedhere\par
\end{samepage}
\end{proof}

\begin{corollary}\label{cor:finiteentailment}
For all $A,B\in\For$,
\begin{enumerate}
\item[(1)] $A\entfde B$ if and only if \eqref{eq:entt} and
\eqref{eq:entone} hold at every point of every model on $\Frm_{2}$;
\item[(2)] $A\entu B$ if and only if \eqref{eq:entt} and
\eqref{eq:entone} hold at every point of every model on $\Frm_{3}$.
\end{enumerate}
Consequently, $A\entn B$ implies both $A\entfde B$ and $A\entu B$.
\end{corollary}

\begin{proof}
Fix either frame and corresponding valuations $V$ and $v$. By
Lemma~\ref{lem:finitevaluations}, the truth sets are ordered by inclusion
in the same way as the finite component values. Hence
\[
\bigl(\Pset{A}{t}\subseteq\Pset{B}{t}\ \text{ and }
\Pset{A}{1}\subseteq\Pset{B}{1}\bigr)
\iff v(A)\leqslant v(B).
\]
These inclusions state precisely that \eqref{eq:entt} and
\eqref{eq:entone} hold at every point of the model.
The bijections of Lemma~\ref{lem:finitevaluations} allow us to quantify
equivalently over all relational or all finite-valued valuations, proving
(1) and (2). The final implication follows because both frames belong to
the class over which $\entn$ is defined.
\end{proof}

\section{The component systems}\label{sec:components}

Each of $\entt$ and $\entone$ is axiomatised by a binary consequence calculus,
and for each the canonical construction goes through. The ontological case is
treated in full; the epistemic one is its mirror image and is stated at the end
of the section.

\subsection[The ontological component calculus]{The calculus $\NDstt$}\label{sec:calculus-t}

The postulates of $\NDstt$ are the axioms and rules of Table~\ref{tab:post-t},
in which $A,B,C$ range over $\For$.

\begin{table}[ht]
\centering
\begin{tabular}{@{}llll@{}}
$(C1)$ & $A\wedge B\dert A$ &
$(CR1)$ & $A\dert B,\ A\dert C\ /\ A\dert B\wedge C$\\
$(C2)$ & $A\wedge B\dert B$ &
$(CR2)$ & $A\dert C,\ B\dert C\ /\ A\vee B\dert C$\\
$(C3)$ & $A\dert A\vee B$ &
$(CR3)$ & $A\dert B,\ B\dert C\ /\ A\dert C$\\
$(C4)$ & $B\dert A\vee B$ &
$(CR4)$ & $A\dert B\ /\ \negt B\dert\negt A$\\
$(C5)$ & $\nego\nego A\dashv\vdash_{t}A$ & &\\
$(C6)$ & $\negt\negt A\dashv\vdash_{t}A$ & &\\
$(C7)$ & $A\dashv\vdash_{t}\nego A$ & &
\end{tabular}
\caption{Postulates of $\NDstt$.}\label{tab:post-t}
\end{table}
\FloatBarrier

The postulates of $\NDstone$ are obtained from Table~\ref{tab:post-t}
by replacing $\dert$ with $\derone$ and interchanging $\negt$ and $\nego$.
For $i\in\{t,1\}$, a derivation in the component calculus is a finite
sequence of formula pairs written $A\vdash_i B$. Each pair is an axiom
instance or follows from earlier pairs by a rule instance. A pair is
derivable if it is the last pair of such a derivation. The notation
$A\dashv\vdash_i B$ stands for both $A\vdash_i B$ and $B\vdash_i A$.

Postulate $(C7)$ is the syntactic counterpart of the transparency conditions
\eqref{eq:cl13}: the epistemic negation leaves the ontological component
unchanged. Applying $(C7)$ to $A$ and $\nego A$ and using $(CR3)$ derives
$(C5)$; we retain it so that both double-negation laws are listed explicitly.
The rule $(CR4)$ is contraposition for $\negt$; its failure for $\entn$
is established in Lemma~\ref{lem:nocontraposition}.

\begin{lemma}\label{lem:soundness-t}
Every axiom of Table~\ref{tab:post-t} is valid for $\entt$, and every rule
preserves validity for $\entt$.
\end{lemma}

\begin{proof}
By Definition~\ref{def:entailment} it suffices to verify the inclusion of
ontological truth sets in every model. By
Corollary~\ref{cor:homomorphism}, conjunction and disjunction are
interpreted by meet and join in $\Lstab_t$. Their defining properties
give $(C1)$--$(C4)$ and $(CR1)$--$(CR2)$; transitivity of inclusion gives
$(CR3)$. The transparency conditions \eqref{eq:cl13} give $(C5)$ and $(C7)$.

For $(C6)$, Corollary~\ref{cor:homomorphism} and the involutivity of
$\Nc{t}$ from Lemma~\ref{lem:starinvolution} give
\[
\Pset{\negt\negt A}{t}
=\Nc{t}\bigl(\Nc{t}(\Pset{A}{t})\bigr)
=\Pset{A}{t}.
\]
For $(CR4)$, assume $\Pset{A}{t}\subseteq\Pset{B}{t}$. By the same
corollary and the antitonicity of $\Nc{t}$,
\[
\Pset{\negt B}{t}
=\Nc{t}(\Pset{B}{t})
\subseteq\Nc{t}(\Pset{A}{t})
=\Pset{\negt A}{t}.
\]
Since the model was arbitrary, every axiom is valid and every rule
preserves validity for $\entt$.
\end{proof}

\begin{theorem}\label{thm:soundness-t}
If $A\dert B$ then $A\entt B$.
\end{theorem}

\begin{proof}
Induction on the length of a derivation, using Lemma~\ref{lem:soundness-t}.
\end{proof}

\subsection{Theories, countertheories and quadruples}\label{sec:theories}

The canonical frame will have points built from disjoint
theory--countertheory pairs, one for each component. The abbreviation $tc$
stands for theory--countertheory.

\begin{definition}\label{def:theory}
Let $i\in\{t,1\}$. A set $T\subseteq\For$ is an
\emph{$\NDst^{i}$-theory} if $A\in T$ and $A\vdash_i B$ imply $B\in T$,
and $A,B\in T$ imply $A\wedge B\in T$. A set $T^{\circ}\subseteq\For$
is an \emph{$\NDst^{i}$-countertheory} if $A\in T^{\circ}$ and
$B\vdash_i A$ imply $B\in T^{\circ}$, and $A,B\in T^{\circ}$ imply
$A\vee B\in T^{\circ}$.
\end{definition}

\begin{definition}\label{def:quadruple}
A \emph{$tc$-quadruple} is a tuple
$x=\langle\Th{x}{t},\Th{x}{1},\Cth{x}{t},\Cth{x}{1}\rangle$ in which,
for each $i\in\{t,1\}$, $\Th{x}{i}$ is an $\NDst^{i}$-theory,
$\Cth{x}{i}$ is an $\NDst^{i}$-countertheory, and
$\Th{x}{i}\cap\Cth{x}{i}=\emptyset$. For $tc$-quadruples $x,y$ put
$x\subseteq y$ if the inclusion holds in each of the four components. A
$tc$-quadruple is \emph{maximal} if it is maximal under $\subseteq$.
\end{definition}

\begin{samepage}
\begin{lemma}\label{lem:extension}
Every $tc$-quadruple is contained in a maximal one.
\end{lemma}

\begin{proof}
Fix a $tc$-quadruple $x$. The set of its extensions is nonempty, since it
contains $x$. For any nonempty chain of extensions, take unions in each
coordinate. The closure conditions are preserved. Each involves only
finitely many formulae, which occur together in one chain member.
Disjointness is preserved as well. A formula in both corresponding unions
would lie in both coordinates of one chain member. Thus every chain has an upper bound, and
Zorn's lemma gives a maximal extension of $x$.
\end{proof}
\end{samepage}

\subsection{The canonical frame}\label{sec:canonical-t}

\begin{definition}\label{def:canonical}
The \emph{canonical frame} is the structure
$\Fc=\langle W,\lt^{c},\rt^{c},\lo^{c},\ro^{c},\Nt^{c},\No^{c}\rangle$ in which
$W$ is the set of all maximal $tc$-quadruples and, for $x,y\in W$,
\begin{equation}\label{eq:canorders}
\begin{aligned}
x\lt^{c}y&\defiff\Th{x}{t}\subseteq\Th{y}{t},
&\qquad
x\rt^{c}y&\defiff\Cth{x}{t}\subseteq\Cth{y}{t},\\
x\lo^{c}y&\defiff\Th{x}{1}\subseteq\Th{y}{1},
&\qquad
x\ro^{c}y&\defiff\Cth{x}{1}\subseteq\Cth{y}{1},
\end{aligned}
\end{equation}
and
\begin{equation}\label{eq:canstars}
\Nt^{c}x=\langle\negt\Cth{x}{t},\Th{x}{1},\negt\Th{x}{t},\Cth{x}{1}\rangle,
\qquad
\No^{c}x=\langle\Th{x}{t},\nego\Cth{x}{1},\Cth{x}{t},\nego\Th{x}{1}\rangle,
\end{equation}
where the operation on sets of formulae is defined by
\[
\neg_iY:=\{A\in\For\mid\neg_iA\in Y\},
\qquad Y\subseteq\For,\quad i\in\{t,1\}.
\]
\end{definition}

\begin{lemma}\label{lem:canwelldef}
$\Fc$ is an $NDst_{oe}$-frame.
\end{lemma}

\begin{proof}
The quadruple of four empty sets has a maximal extension by
Lemma~\ref{lem:extension}, so $W\ne\emptyset$. The four relations
\eqref{eq:canorders} are quasi-orders, being inclusions.

In each component, the lattice postulates and contraposition give one
direction of each De Morgan law. Substituting $\neg_i A$ and $\neg_i B$
in these sequents and using the double-negation and lattice rules gives
\[
\neg_i(\neg_i A\vee\neg_i B)\vdash_i A\wedge B,
\qquad
A\vee B\vdash_i\neg_i(\neg_i A\wedge\neg_i B).
\]
Applying contraposition and eliminating double negations gives the
remaining directions. Thus
\[
\begin{aligned}
\neg_i(A\vee B)&\dashv\vdash_i\neg_iA\wedge\neg_iB,\\
\neg_i(A\wedge B)&\dashv\vdash_i\neg_iA\vee\neg_iB.
\end{aligned}
\]
If $T$ is an $\NDst^{i}$-theory, contraposition gives downward closure of
$\neg_iT$. The first De Morgan law gives its closure under $\vee$.
Thus $\neg_iT$ is an $\NDst^{i}$-countertheory. The dual argument shows
that $\neg_iT^{\circ}$ is a theory whenever $T^{\circ}$ is a countertheory.

Double negation gives $\neg_i\neg_iY=Y$ for every theory or countertheory $Y$.
The operation preserves inclusion and disjointness. Hence the formulae in
\eqref{eq:canstars} define inclusion-preserving involutions on all
$tc$-quadruples. They preserve maximality and therefore map $W$ to itself.
This also proves \eqref{eq:frinv}.

By involutivity, the operation also reflects inclusion between such sets.
Together with \eqref{eq:canorders} and \eqref{eq:canstars}, this gives
\[
x\ordl{i}^{c}y\iff N_i^{c}x\ordr{i}^{c}N_i^{c}y,
\]
which is \eqref{eq:frswap}. Finally, $\Nt^{c}$ and $\No^{c}$ commute because
they act on different pairs of coordinates.
\end{proof}

The next lemma gives the conjugacy conditions for the canonical valuation.

\begin{lemma}\label{lem:cansat}
For every $x\in W$ and every $A\in\For$,
\begin{equation}\label{eq:cansat}
\begin{aligned}
A\in\Th{x}{t}&\iff\forall y\,(x\lt^{c}y\Rightarrow A\notin\Cth{y}{t}),
&\qquad
A\in\Cth{x}{t}&\iff\forall y\,(x\rt^{c}y\Rightarrow A\notin\Th{y}{t}),\\
A\in\Th{x}{1}&\iff\forall y\,(x\lo^{c}y\Rightarrow A\notin\Cth{y}{1}),
&\qquad
A\in\Cth{x}{1}&\iff\forall y\,(x\ro^{c}y\Rightarrow A\notin\Th{y}{1}).
\end{aligned}
\end{equation}
\end{lemma}

\begin{proof}
We prove the first equivalence. If $A\in\Th{x}{t}$ and $x\lt^{c}y$, then
$A\in\Th{y}{t}$, so $A\notin\Cth{y}{t}$ by disjointness.

Conversely, suppose $A\notin\Th{x}{t}$ and put
$\Gamma:=\{B\in\For\mid B\dert A\}$. By $(C6)$ and $(CR3)$,
$A\dert A$, so $A\in\Gamma$. Rule $(CR3)$ gives downward closure of
$\Gamma$, and $(CR2)$ gives closure under $\vee$. Thus $\Gamma$ is a
countertheory. It is disjoint from $\Th{x}{t}$ by the closure of
$\Th{x}{t}$ under derivability. Apply Lemma~\ref{lem:extension} to
$\langle\Th{x}{t},\Th{x}{1},\Gamma,\Cth{x}{1}\rangle$. The resulting
maximal quadruple $y$ satisfies $x\lt^{c}y$ and $A\in\Cth{y}{t}$.
The other equivalences follow by exchanging theories with countertheories
and interchanging the two components.
\end{proof}

\begin{corollary}\label{cor:canvaluation}
Setting $\Pset{p}{t}=\{x\in W\mid p\in\Th{x}{t}\}$,
$\Pset{p}{f}=\{x\in W\mid p\in\Cth{x}{t}\}$,
$\Pset{p}{1}=\{x\in W\mid p\in\Th{x}{1}\}$ and
$\Pset{p}{0}=\{x\in W\mid p\in\Cth{x}{1}\}$ defines a valuation on $\Fc$.
\end{corollary}

\begin{proof}
By \eqref{eq:cansat} and Definition~\ref{def:galoismaps},
$\rmap{t}(\Pset{p}{t})=\Pset{p}{f}$ and $\lmap{t}(\Pset{p}{f})=\Pset{p}{t}$,
which is $t$-conjugacy; the epistemic case is the same.
\end{proof}

\subsection{Completeness}\label{sec:completeness-t}

\begin{lemma}\label{lem:truth}
In the canonical model of Corollary~\ref{cor:canvaluation}, for every
$A\in\For$ and every $x\in W$,
\[
x\forc{t}A\iff A\in\Th{x}{t},
\qquad
x\forc{f}A\iff A\in\Cth{x}{t},
\]
and correspondingly for $\forc{1}$ and $\forc{0}$.
\end{lemma}

\begin{proof}
We prove the four equivalences simultaneously by structural induction.
The variable case is Corollary~\ref{cor:canvaluation}. At each step we
first establish the truth equivalences at every point.
Lemma~\ref{lem:conjugacypreserved} and \eqref{eq:cansat} then give the
falsity equivalences: both falsity sets are obtained by applying $\rmap{i}$
to the same truth set.

We give the ontological truth cases. For $A\wedge B$, closure of theories
under $\wedge$ together with $(C1)$ and $(C2)$ gives
$A\wedge B\in\Th{x}{t}$ if and only if $A,B\in\Th{x}{t}$.
The induction hypotheses and \eqref{eq:cl1} give the result.

By downward closure and closure under $\vee$, a countertheory contains
$A\vee B$ if and only if it contains both $A$ and $B$. Using the induction
hypotheses, \eqref{eq:cl5} therefore says that no $\lt^{c}$-successor of $x$
has $A\vee B$ in its countertheory. By \eqref{eq:cansat}, this is equivalent
to $A\vee B\in\Th{x}{t}$.

For $\negt A$, \eqref{eq:cl9}, the induction hypothesis and
\eqref{eq:canstars} give $x\forc{t}\negt A$ if and only if
$A\in\negt\Th{x}{t}$, which means $\negt A\in\Th{x}{t}$ by definition.
For $\nego A$, the transparency conditions \eqref{eq:cl13} and $(C7)$ give
the result. The epistemic truth cases follow by interchanging the two
components and the two negations.
\end{proof}

\begin{theorem}\label{thm:completeness-t}
$A\dert B$ if and only if $A\entt B$.
\end{theorem}

\begin{proof}
Soundness is Theorem~\ref{thm:soundness-t}. For the converse, suppose
$A\not\dert B$ and put
\[
T:=\{C\in\For\mid A\dert C\},
\qquad I:=\{C\in\For\mid C\dert B\}.
\]
Then $T$ is a theory and $I$ is a countertheory by $(CR1)$--$(CR3)$.
They are disjoint by $(CR3)$, and reflexivity gives $A\in T$ and $B\in I$.
By Lemma~\ref{lem:extension}, the $tc$-quadruple
$\langle T,\emptyset,I,\emptyset\rangle$ has a maximal extension $x$.
Thus $A\in\Th{x}{t}$ and $B\in\Cth{x}{t}$.
Disjointness and Lemma~\ref{lem:truth} give $x\forc{t}A$ and
$x\nforc{t}B$. Hence $A\not\entt B$ by Definition~\ref{def:entailment}.
\end{proof}

\begin{remark}\label{rem:mirror}
The arguments of Sections~\ref{sec:calculus-t}--\ref{sec:completeness-t}
carry over by interchanging the ontological and epistemic components.
In particular,
$A\derone B$ if and only if $A\entone B$.
\end{remark}

\section[The system NDst\_oe]{The system $\NDst$}\label{sec:system}

This section presents a calculus for $\entn=\entt\cap\entone$.
We prove its soundness and its inclusion in both component calculi.
Section~\ref{sec:algebra} establishes the converse inclusion and completeness.

\subsection{Structural properties}\label{sec:joint-properties}

Neither negation supports unrestricted contraposition for $\entn$.

\begin{lemma}\label{lem:nocontraposition}
There are $A,B\in\For$ with $A\entn B$ and $\negt B\not\entn\negt A$.
\end{lemma}

\begin{proof}
Take $A:=p$ and $B:=p\vee q$. By Corollary~\ref{cor:homomorphism} and
Definition~\ref{def:entailment}, $p\entn p\vee q$.
In $\Four$, take $v(p)=\langle f,0\rangle$ and
$v(q)=\langle t,1\rangle$. By \eqref{eq:fourops},
\[
v(\negt(p\vee q))=\langle f,1\rangle
\not\leqslant\langle t,0\rangle=v(\negt p).
\]
Thus $\negt(p\vee q)\not\entfde\negt p$.
Corollary~\ref{cor:finiteentailment} gives
$\negt(p\vee q)\not\entn\negt p$.
\end{proof}

By the transparency condition \eqref{eq:cl14}, $\negt$ leaves epistemic truth
sets unchanged. Thus a strict inclusion
$\Pset{A}{1}\subsetneq\Pset{B}{1}$ in some model refutes
$\negt B\entn\negt A$.
The argument for $\nego$ in the ontological component is symmetric.

The relational semantics also differs from the finite-valued semantics
of Section~\ref{sec:finite}: distributivity need not hold.

\begin{lemma}\label{lem:nondistributive}
The schema
\begin{equation}\label{eq:distribution}
A\wedge(B\vee C)\ \entn\ (A\wedge B)\vee(A\wedge C)
\end{equation}
is not valid.
\end{lemma}

\begin{proof}
Take $X=\{x,y,z\}$. In each component, the only nonreflexive comparisons are
\[
x\ordl{i}z,\qquad x\ordr{i}y,
\qquad i\in\{t,1\}.
\]
Let both $\Nt$ and $\No$ fix $x$ and interchange $y,z$.
They are commuting involutions that exchange the two quasi-orders.
Thus the structure is an $NDst_{oe}$-frame.

The $\lmap{i}$-stable subsets of $X$ form the lattice
\[
\Lstab_i=\{\emptyset,\{z\},\{x,z\},\{y\},X\}.
\]
Indeed, the only other $\ordl{i}$-increasing subset is $\{y,z\}$.
The empty set and $X$ are stable. By Definition~\ref{def:galoismaps},
$\rmap{i}$ sends $\{z\}$ to $\{x,y\}$, $\{x,z\}$ to $\{y\}$, and
$\{y\}$ to $\{z\}$. In each case $\lmap{i}$ returns the original set.
Finally, $\rmap{i}(\{y,z\})=\emptyset$, so the closure of $\{y,z\}$ is $X$.

Under inclusion, this is the pentagon $N_5$.
Put $a:=\{z\}$, $b:=\{x,z\}$ and $c:=\{y\}$.
Its lattice operations satisfy
\[
b\wedge(c\vee a)=b\not\subseteq a=(b\wedge c)\vee(b\wedge a).
\]

Take distinct propositional variables $p,q,r$.
In both components, assign them the truth sets $b,c,a$, respectively.
Lemma~\ref{lem:conjugate} gives the corresponding conjugate pairs.
By Corollary~\ref{cor:homomorphism}, $p\wedge(q\vee r)$ and
$(p\wedge q)\vee(p\wedge r)$ have ontological truth sets $b$ and $a$.
Since $x\in b\setminus a$, this instance of \eqref{eq:distribution}
fails at $x$.
\end{proof}

The construction in Section~\ref{sec:canonical-t} relies on contraposition.
Section~\ref{sec:algebra} uses two quotients of the Lindenbaum algebra,
one for each component.

\subsection{The calculus}\label{sec:calculus-n}

The postulates of $\NDst$ are listed in Table~\ref{tab:post-n} and
\eqref{eq:rg}. Here $A,B,C$ range over $\For$.
Derivations are defined in the same way as in the component calculi,
with formula pairs written $A\dern B$. The notation $A\dashv\vdash_n B$
stands for both $A\dern B$ and $B\dern A$. In $(NR6)$ and $(NR7)$, both
directions of the premiss are required, and each direction of the
conclusion may be inferred separately.

\begin{table}[H]
\centering
\begin{tabular}{@{}llll@{}}
$(N1)$ & $A\wedge B\dern A$ &
$(NR1)$ & $A\dern B,\ A\dern C\ /\ A\dern B\wedge C$\\
$(N2)$ & $A\wedge B\dern B$ &
$(NR2)$ & $A\dern C,\ B\dern C\ /\ A\vee B\dern C$\\
$(N3)$ & $A\dern A\vee B$ &
$(NR3)$ & $A\dern B,\ B\dern C\ /\ A\dern C$\\
$(N4)$ & $B\dern A\vee B$ &
$(NR4)$ & $\negt A\dern\nego B\ /\ \negt B\dern\nego A$\\
$(N5)$ & $\nego\nego A\dashv\vdash_{n}A$ &
$(NR5)$ & $\nego A\dern\negt B\ /\ \nego B\dern\negt A$\\
$(N6)$ & $\negt\negt A\dashv\vdash_{n}A$ &
$(NR6)$ & $A\dashv\vdash_{n}B\ /\ \negt A\dashv\vdash_{n}\negt B$\\
$(N7)$ & $\negt A\wedge B\dern\negt(A\wedge B)$ &
$(NR7)$ & $A\dashv\vdash_{n}B\ /\ \nego A\dashv\vdash_{n}\nego B$\\
$(N8)$ & $\nego A\wedge B\dern\nego(A\wedge B)$ & &\\
$(N9)$ & $\negt(A\wedge B)\dern\negt A\vee\negt B$ & &\\
$(N10)$ & $\nego(A\wedge B)\dern\nego A\vee\nego B$ & &
\end{tabular}
\caption{Axioms $(N1)$--$(N10)$ and rules $(NR1)$--$(NR7)$ of
$\NDst$.}\label{tab:post-n}
\end{table}

\noindent Each of $(NR8)$ and $(NR9)$ abbreviates two rules with the same premisses.
\begin{equation}\label{eq:rg}
\begin{gathered}
(NR8)\ \
\frac{A\dern B,\quad\negt A\dern\negt B}
{\negt(A\wedge C)\dern\negt(B\wedge C),\quad\negt(A\vee C)\dern\negt(B\vee C)}
\\[8pt]
(NR9)\ \
\frac{A\dern B,\quad\nego A\dern\nego B}
{\nego(A\wedge C)\dern\nego(B\wedge C),\quad\nego(A\vee C)\dern\nego(B\vee C)}
\end{gathered}
\end{equation}

Unrestricted contraposition $(CR4)$ is replaced by the cross rules
$(NR4)$ and $(NR5)$. The transparency postulate $(C7)$ is also absent:
neither negation is transparent for $\entn$.
The rules $(NR8)$ and $(NR9)$ will ensure compatibility of the component
congruences with the lattice operations; see
Lemma~\ref{lem:component-congruences}.

\subsection{Soundness}\label{sec:soundness-n}

We verify soundness separately for the two component entailment relations
$\entt$ and $\entone$.
For $i\in\{t,1\}$, recall that $\Pset{A}{i}=\{x\mid x\forc{i}A\}$.
Thus, by Definition~\ref{def:entailment}, $A\models^{i}_{oe}B$ means that
$\Pset{A}{i}\subseteq\Pset{B}{i}$ in every model. By
Corollary~\ref{cor:homomorphism}, in each component one negation is
transparent and the other acts as an antitone involution of the lattice
of stable sets.

\begin{lemma}\label{lem:soundness-n}
Every axiom of Table~\ref{tab:post-n} is valid for $\entt$ and for $\entone$,
and every rule of Table~\ref{tab:post-n} and \eqref{eq:rg} preserves validity
for each of them.
\end{lemma}

\begin{proof}
Fix an arbitrary model. We verify the ontological component.
The epistemic case follows by exchanging the two components and the two
negations.

Postulates $(N1)$--$(N4)$ and the rules $(NR1)$--$(NR3)$ hold because, by
Corollary~\ref{cor:homomorphism}, $\wedge$ and $\vee$ are computed as the meet
and the join of $\Lstab_{t}$, and inclusion is a partial order. $(N5)$ follows from
Lemma~\ref{lem:negcomponents}(1) applied twice, and $(N6)$ from the
involutivity of $\Nc{t}$.

$(N7)$: from $\Pset{A\wedge B}{t}\subseteq\Pset{A}{t}$ and the antitonicity of
$\Nc{t}$ we obtain $\Pset{\negt A}{t}\subseteq\Pset{\negt(A\wedge B)}{t}$, and
$\Pset{\negt A\wedge B}{t}\subseteq\Pset{\negt A}{t}$. $(N8)$ is an equality by
Lemma~\ref{lem:negcomponents}(1). $(N9)$ follows from the first identity of
\eqref{eq:demorgan}, and $(N10)$ reduces by
Lemma~\ref{lem:negcomponents}(1) to $\Pset{A\wedge B}{t}\subseteq
\Pset{A\vee B}{t}$.

$(NR4)$: by Lemma~\ref{lem:negcomponents}(1) the premiss is
$\Pset{\negt A}{t}\subseteq\Pset{B}{t}$; applying $\Nc{t}$ and using its
antitonicity and involutivity gives
$\Pset{\negt B}{t}\subseteq\Pset{A}{t}=\Pset{\nego A}{t}$. $(NR5)$: the premiss
is $\Pset{A}{t}\subseteq\Pset{\negt B}{t}$ by the same condition, and $\Nc{t}$
turns it into $\Pset{B}{t}\subseteq\Pset{\negt A}{t}$, which is the conclusion.

$(NR6)$: the premiss gives $\Pset{A}{t}=\Pset{B}{t}$.
Applying $\Nc{t}$ preserves this equality. $(NR7)$ follows from
Lemma~\ref{lem:negcomponents}(1).

$(NR8)$: the second premiss reads
$\Nc{t}(\Pset{A}{t})\subseteq\Nc{t}(\Pset{B}{t})$.
Antitonicity and involutivity of $\Nc{t}$ yield
$\Pset{B}{t}\subseteq\Pset{A}{t}$.
Together with the first premiss, this gives $\Pset{A}{t}=\Pset{B}{t}$.
Hence $\Pset{A\wedge C}{t}=\Pset{B\wedge C}{t}$ and
$\Pset{A\vee C}{t}=\Pset{B\vee C}{t}$.
Applying $\Nc{t}$ gives both conclusions.

$(NR9)$: by Lemma~\ref{lem:negcomponents}(1), the second premiss coincides
with the first. The conclusions reduce to
$\Pset{A\wedge C}{t}\subseteq\Pset{B\wedge C}{t}$ and
$\Pset{A\vee C}{t}\subseteq\Pset{B\vee C}{t}$.
These follow from the first premiss by monotonicity of the lattice operations.
\end{proof}

\begin{theorem}\label{thm:soundness-n}
If $A\dern B$ then $A\entn B$.
\end{theorem}

\begin{proof}
Induction on the length of a derivation gives $A\dern B\Rightarrow A\entt B$
and $A\dern B\Rightarrow A\entone B$ by Lemma~\ref{lem:soundness-n}; now apply
\eqref{eq:semdecomp}.
\end{proof}

\subsection{Inclusion in the component calculi}\label{sec:embedding}

We prove inclusion in both component calculi by syntactic means.
We first establish some basic properties of $\NDstt$.

\begin{lemma}\label{lem:derivedt}
In $\NDstt$, if $A\dert B$, then $A\wedge C\dert B\wedge C$ and
$A\vee C\dert B\vee C$.
If $A\dashv\vdash_{t}B$, then $\negt A\dashv\vdash_{t}\negt B$ and
$\nego A\dashv\vdash_{t}\nego B$.
Moreover, $\nego A\dashv\vdash_{t}A$ and
$\negt(A\wedge B)\dashv\vdash_{t}\negt A\vee\negt B$.
\end{lemma}

\begin{proof}
From $A\dert B$ and $(C1)$, $(C2)$, $(CR3)$ we get $A\wedge C\dert B$ and
$A\wedge C\dert C$, whence $A\wedge C\dert B\wedge C$ by $(CR1)$.
The case of $\vee$ is dual, using $(C3)$, $(C4)$, $(CR2)$ and $(CR3)$.
If also $B\dert A$, the converse derivations follow in the same way.
For $\negt$ apply $(CR4)$ to each of $A\dert B$ and
$B\dert A$; for $\nego$ use $(C7)$ twice together with $(CR3)$. The statement
$\nego A\dashv\vdash_{t}A$ is $(C7)$.

For the De Morgan law, $(C1)$, $(C2)$ and $(CR4)$ give
$\negt A\dert\negt(A\wedge B)$ and $\negt B\dert\negt(A\wedge B)$, whence
$\negt A\vee\negt B\dert\negt(A\wedge B)$ by $(CR2)$. Conversely, $(C3)$,
$(C4)$ and $(CR4)$ give $\negt(\negt A\vee\negt B)\dert\negt\negt A$ and
$\negt(\negt A\vee\negt B)\dert\negt\negt B$, hence
$\negt(\negt A\vee\negt B)\dert A\wedge B$ by $(C6)$ and $(CR1)$.
A further application of $(CR4)$ and $(C6)$ yields
$\negt(A\wedge B)\dert\negt A\vee\negt B$.
\end{proof}

\begin{lemma}\label{lem:postulates-in-components}
Every axiom of Table~\ref{tab:post-n} is derivable in $\NDstt$ and in
$\NDstone$, and every rule of Table~\ref{tab:post-n} and \eqref{eq:rg} is
admissible in each of them.
\end{lemma}

\begin{proof}
We argue for $\NDstt$.
Interchanging the two negations preserves the postulates of $\NDst$.
By the definition of $\NDstone$ in Section~\ref{sec:calculus-t}, the same
exchange gives the proof for $\NDstone$.

The axioms $(N1)$--$(N4)$ are $(C1)$--$(C4)$.
Axioms $(N5)$ and $(N6)$ are $(C5)$ and $(C6)$.
By $(C7)$ and Lemma~\ref{lem:derivedt}, $(N8)$ reduces to
$A\wedge B\dert A\wedge B$, which holds by reflexivity.
Similarly, $(N10)$ reduces to $A\wedge B\dert A\vee B$.
This follows from $(C1)$, $(C3)$ and $(CR3)$.
Axiom $(N9)$ follows from Lemma~\ref{lem:derivedt}.
For $(N7)$, applying $(CR4)$ to $(C1)$ gives
$\negt A\dert\negt(A\wedge B)$.
Compose this with $\negt A\wedge B\dert\negt A$, using $(C1)$ and $(CR3)$.

The rules $(NR1)$--$(NR3)$ are $(CR1)$--$(CR3)$.
Rules $(NR6)$ and $(NR7)$ follow from Lemma~\ref{lem:derivedt}.

For $(NR4)$, the premiss $\negt A\dert\nego B$ is equivalent to
$\negt A\dert B$ by $(C7)$ and $(CR3)$.
Applying $(CR4)$ gives $\negt B\dert\negt\negt A$.
By $(C6)$, $(C7)$ and $(CR3)$, this yields $\negt B\dert\nego A$.

For $(NR5)$, the premiss is equivalent to $A\dert\negt B$ by $(C7)$ and $(CR3)$.
Applying $(CR4)$ and using $(C6)$ gives $B\dert\negt A$.
By $(C7)$ and $(CR3)$, this yields $\nego B\dert\negt A$.

For $(NR8)$, applying $(CR4)$ to the second premiss and using $(C6)$ gives
$B\dert A$. Thus $A\dashv\vdash_{t}B$.
By Lemma~\ref{lem:derivedt}, the
formulae $A\wedge C$ and $B\wedge C$, and likewise $A\vee C$ and $B\vee C$, are
interderivable.
Preservation of interderivability under $\negt$ gives both conclusions.

For $(NR9)$, the second premiss is equivalent to the first by $(C7)$.
The conclusions reduce to $A\wedge C\dert B\wedge C$ and
$A\vee C\dert B\vee C$ by $(C7)$ and $(CR3)$.
Both follow from the first premiss by Lemma~\ref{lem:derivedt}.
\end{proof}

\begin{corollary}\label{cor:embedding}
$\dern\ \subseteq\ \dert\cap\derone$.
\end{corollary}

\begin{proof}
Induction on the length of a derivation in $\NDst$, using
Lemma~\ref{lem:postulates-in-components} at each step.
\end{proof}

The converse inclusion is established algebraically in
Section~\ref{sec:algebra}. The syntactic decomposition is independent of
the component completeness theorem.

\section{Algebraic decomposition and completeness}\label{sec:algebra}

We construct two quotients of the Lindenbaum algebra of $\NDst$,
one for each component calculus.
Together, the quotient maps reflect the order of the original algebra.

\subsection{The Lindenbaum algebra}\label{sec:lindenbaum}

\begin{definition}\label{def:joint-lindenbaum}
For $A,B\in\For$, put $A\equiv_n B$ if $A\dashv\vdash_n B$.
Let $L:=\For/{\equiv_n}$ and write $[A]$ for the class of $A$.
Set
\[
[A]\leqslant[B]\defiff A\dern B,
\qquad
[A]\wedge[B]:=[A\wedge B],
\qquad
[A]\vee[B]:=[A\vee B].
\]
We also use $\negt,\nego$ for the induced operations:
$\negt[A]:=[\negt A]$ and $\nego[A]:=[\nego A]$.
\end{definition}

\begin{lemma}\label{lem:lindenbaum-structure}
The relation $\equiv_n$ is a congruence on the formula algebra.
The order and operations above are well-defined, and $L$ is a lattice.
The operations $\negt,\nego$ are commuting involutions, and their composition
$k:=\negt\circ\nego$ is an antitone involution of $L$.
\end{lemma}

\begin{proof}
Reflexivity of $\dern$ follows from $(N6)$ and $(NR3)$;
transitivity follows from $(NR3)$. Hence $\equiv_n$ is an equivalence relation.
Suppose $A\equiv_n A'$ and $B\equiv_n B'$.
By $(N1)$, $(N2)$ and $(NR3)$, $A\wedge B\dern A'$ and
$A\wedge B\dern B'$. Rule $(NR1)$ gives $A\wedge B\dern A'\wedge B'$.
The reverse derivation is analogous.
The case of $\vee$ is dual, using $(N3)$, $(N4)$, $(NR2)$ and $(NR3)$.
Compatibility with the two negations follows from $(NR6)$ and $(NR7)$.
Thus $\equiv_n$ is a congruence, and the operations on classes are well-defined.

If $A\equiv_n A'$ and $B\equiv_n B'$, then $(NR3)$ gives
$A\dern B$ if and only if $A'\dern B'$.
Thus the order on classes is well-defined.
Reflexivity and transitivity follow from those of $\dern$;
antisymmetry follows from the definition of $\equiv_n$.
By $(N1)$ and $(N2)$, $[A\wedge B]$ is a lower bound of $[A]$ and $[B]$.
If $[C]\leqslant[A]$ and $[C]\leqslant[B]$, then $(NR1)$ gives
$C\dern A\wedge B$, so $[C]\leqslant[A\wedge B]$.
The dual argument, using $(N3)$, $(N4)$ and $(NR2)$, shows that $[A\vee B]$
is their least upper bound. Hence $L$ is a lattice.

Axioms $(N5)$ and $(N6)$ give $\negt\negt a=\nego\nego a=a$ for every $a\in L$.
Together with $(NR3)$, they also give
$\negt\negt A\dern\nego\nego A$ and
$\nego\nego A\dern\negt\negt A$.
Applying $(NR4)$ to the first pair and $(NR5)$ to the second yields
\[
\negt\nego A\dern\nego\negt A,
\qquad
\nego\negt A\dern\negt\nego A.
\]
Thus $\negt\circ\nego=\nego\circ\negt$.

Suppose $A\dern B$. By $(N5)$, $(N6)$ and $(NR3)$,
$\negt\negt A\dern\nego\nego B$.
Applying $(NR4)$ gives $\negt\nego B\dern\nego\negt A$.
By commutation, $\negt\nego B\dern\negt\nego A$.
Hence $k$ is antitone. Finally, for every $a\in L$,
\[
k^2a=\negt\nego\negt\nego a=\negt\negt\nego\nego a=a.\qedhere
\]
\end{proof}

\subsection{Component congruences}\label{sec:component-congruences}

Fix $i\in\{t,1\}$ and let $a=[A]$ and $b=[B]$. If $a\leqslant b$
and $\neg_i a\leqslant\neg_i b$, then Corollary~\ref{cor:embedding} gives
$A\vdash_i B$ and $\neg_i A\vdash_i\neg_i B$. The latter gives
$B\vdash_i A$ by contraposition and double negation. We write $aP_i b$
when both inequalities hold. The relation $\theta_i$ relates $a$ and $b$
if and only if they have a common $P_i$-upper bound.

\begin{definition}\label{def:component-congruences}
For $a,b\in L$ and $i\in\{t,1\}$, put
\[
aP_i b\defiff a\leqslant b\ \text{ and }\ \neg_i a\leqslant\neg_i b,
\]
and
\[
a\theta_i b\defiff
\exists c\in L\,(aP_i c\ \text{ and }\ bP_i c).
\]
\end{definition}

\begin{lemma}\label{lem:component-congruences}
For each $i\in\{t,1\}$, the relation $\theta_i$ is a congruence of
$\langle L,\wedge,\vee,\negt,\nego\rangle$. Moreover, for all $a,b\in L$,
\[
a\theta_i b\iff\neg_i(a\wedge b)\leqslant\neg_i(a\vee b).
\]
\end{lemma}

\begin{proof}
Fix $i\in\{t,1\}$. The relation $P_i$ is reflexive and transitive.
For every $d\in L$,
\begin{equation}\label{eq:Pi-lattice-compatible}
aP_i b\Longrightarrow
(a\wedge d)P_i(b\wedge d)
\quad\text{and}\quad
(a\vee d)P_i(b\vee d).
\end{equation}
The inequalities under $\neg_i$ in the two conclusions follow from $(NR8)$
for $i=t$ and $(NR9)$ for $i=1$.
The other inequalities follow from monotonicity of the lattice operations.

The relation $\theta_i$ is reflexive (take $c=a$) and symmetric.
For transitivity, suppose $a\theta_i b$ and $b\theta_i c$.
Choose $d,e$ such that
$aP_i d$, $bP_i d$, $bP_i e$ and $cP_i e$.
Since $b\leqslant d$, \eqref{eq:Pi-lattice-compatible} gives
\[
d=b\vee d\;P_i\;e\vee d=d\vee e.
\]
Interchanging $d,e$ gives $eP_i(d\vee e)$.
Transitivity of $P_i$ now yields $aP_i(d\vee e)$ and $cP_i(d\vee e)$.
Thus $a\theta_i c$.

Let $a\theta_i b$, witnessed by $aP_i c$ and $bP_i c$.
For every $d\in L$, \eqref{eq:Pi-lattice-compatible} shows that
$c\wedge d$ and $c\vee d$ witness
\[
(a\wedge d)\theta_i(b\wedge d),
\qquad
(a\vee d)\theta_i(b\vee d),
\]
respectively. Commutativity and transitivity then give compatibility
in both arguments.

Involutivity gives $\neg_i aP_i\neg_i c$: its two defining inequalities are
\[
\neg_i a\leqslant\neg_i c,
\qquad
\neg_i\neg_i a=a\leqslant c=\neg_i\neg_i c.
\]
The same holds for $b$, so $\neg_i a\theta_i\neg_i b$.
For the other negation, let $\{i,j\}=\{t,1\}$.
By Lemma~\ref{lem:lindenbaum-structure}, the negations commute and
$k=\negt\circ\nego$ is antitone. Hence $aP_i c$ gives
\[
\neg_j c=k\neg_i c\leqslant k\neg_i a=\neg_j a,
\qquad
\neg_i\neg_j c=kc\leqslant ka=\neg_i\neg_j a.
\]
Thus $\neg_j c$ is a common $P_i$-lower bound of $\neg_j a$ and $\neg_j b$.
To obtain a common $P_i$-upper bound, put $u:=\neg_j a\vee\neg_j b$.
Indeed, \eqref{eq:Pi-lattice-compatible} gives
\[
\neg_j a=\neg_j c\vee\neg_j a
\;P_i\;
\neg_j b\vee\neg_j a=u.
\]
Interchanging $a,b$ gives $\neg_j bP_i u$, so
$\neg_j a\theta_i\neg_j b$.
This proves that $\theta_i$ is a congruence.

For the stated criterion, put $d:=a\wedge b$ and $e:=a\vee b$.
Suppose $a\theta_i b$, and choose $c$ with $aP_i c$ and $bP_i c$.
By \eqref{eq:Pi-lattice-compatible},
\[
d=a\wedge b\;P_i\;c\wedge b=b,
\qquad
d=b\wedge a\;P_i\;c\wedge a=a.
\]
Applying the join implication in \eqref{eq:Pi-lattice-compatible} to $dP_i b$
gives $a=d\vee a\;P_i\;b\vee a=e$.
Together with $dP_i a$, this yields $dP_i e$, hence
$\neg_i d\leqslant\neg_i e$.
Conversely, this inequality and $d\leqslant e$ give $dP_i e$.
Joining with $a$ and $b$, respectively, gives $aP_i e$ and $bP_i e$.
Thus $e$ witnesses $a\theta_i b$.
\end{proof}

The relation $\theta_i$ is the least equivalence relation containing $P_i$.
It contains $P_i$ by taking $c=b$ in Definition~\ref{def:component-congruences}.
Any equivalence relation containing $P_i$ must also relate $a,b$ whenever
$aP_i c$ and $bP_i c$.

\subsection{Component quotients}\label{sec:component-quotients}

For $i\in\{t,1\}$, put $L_i:=L/\theta_i$ and let $q_i:L\to L_i$ be
the quotient map.
Each quotient carries the operations induced from $L$.
Its lattice order is given by $u\leqslant v\defiff u\wedge v=u$ for
$u,v\in L_i$.

\begin{lemma}\label{lem:component-quotients}
On $L_t$, the operation $\nego$ is the identity and $\negt$ is an antitone
involution. On $L_1$, the roles of $\negt$ and $\nego$ are interchanged.
Consequently, for all $A,B\in\For$ and each $i\in\{t,1\}$,
\[
A\vdash_i B\Longrightarrow q_i([A])\leqslant q_i([B]).
\]
\end{lemma}

\begin{proof}
We prove the assertions for $i=t$.
The induced negations are involutive by Lemma~\ref{lem:lindenbaum-structure}.
Write $k:=\negt\circ\nego$.
For $a\in L$, put $d:=a\wedge\nego a$ and $e:=a\vee\nego a$.
Two instances of $(N8)$ give
\[
\begin{aligned}
\nego a&=\nego a\wedge\nego a
\leqslant\nego(a\wedge\nego a)=\nego d,\\
a&=\nego\nego a\wedge a
\leqslant\nego(\nego a\wedge a)=\nego d.
\end{aligned}
\]
Hence $e\leqslant\nego d$.
Axiom $(N10)$ gives $\nego d\leqslant\nego a\vee\nego\nego a=e$.
Thus $\nego d=e$ and, by involutivity, $\nego e=d$.
Since $d\leqslant e$ and $k$ is antitone,
\[
\negt d=ke\leqslant kd=\negt e.
\]
The criterion in Lemma~\ref{lem:component-congruences} now gives
$a\theta_t\nego a$.
Consequently, $\nego q_t(a)=q_t(\nego a)=q_t(a)$.
Surjectivity of $q_t$ shows that $\nego$ is the identity on $L_t$.

For every $x\in L$, we therefore have
\[
q_t(kx)=\negt\nego q_t(x)=\negt q_t(x).
\]
As an antitone involution, $k$ satisfies $k(a\wedge b)=ka\vee kb$.
Applying $q_t$ to this equality and using its surjectivity gives
\[
\negt(u\wedge v)=\negt u\vee\negt v
\qquad(u,v\in L_t).
\]
If $u\leqslant v$, then $u\wedge v=u$, so
$\negt u\vee\negt v=\negt u$.
Thus $\negt v\leqslant\negt u$, proving antitonicity on $L_t$.

For the final assertion, put $w(A):=q_t([A])$.
The map $w$ preserves all operations as a composition of quotient maps.
We prove $w(A)\leqslant w(B)$ by induction on the length of a derivation
of $A\dert B$.
Axioms $(C1)$--$(C4)$ follow from the defining properties of $\wedge$ and $\vee$.
Axioms $(C5)$, $(C6)$ follow from involutivity, and $(C7)$ from the identity
$\nego u=u$ on $L_t$.
For $(CR1)$ and $(CR2)$, apply the induction hypothesis and the properties of
greatest lower bounds and least upper bounds, respectively.
Rule $(CR3)$ follows from the induction hypothesis and transitivity.
For $(CR4)$, the induction hypothesis gives $w(A)\leqslant w(B)$, hence
\[
w(\negt B)=\negt w(B)\leqslant\negt w(A)=w(\negt A).
\]
This completes the induction.
Interchanging the indices $t,1$ gives the assertions for $i=1$.
\end{proof}

The next lemma uses the fact that quotient maps give a subdirect
embedding when their kernels intersect in the
diagonal~\cite[Chapter~II, Lemma~8.2]{BurrisSankappanavar}.
The product representations for interlaced
pre-bilattices~\cite[Proposition~3.8]{BouRivieccio} and interlaced
trilattices~\cite[Section~3]{RivieccioTrilattices} likewise use two
congruences.

\begin{lemma}\label{lem:component-separation}
The two congruences satisfy
\[
\theta_t\cap\theta_1=\Delta_L,
\qquad \Delta_L:=\{(a,a)\mid a\in L\}.
\]
Consequently, the map
\[
q:L\longrightarrow L_t\times L_1,
\qquad q(a):=\bigl(q_t(a),q_1(a)\bigr),
\]
is a subdirect embedding and reflects order.
\end{lemma}

\begin{proof}
Suppose $a\theta_t b$ and $a\theta_1 b$.
Put $d:=a\wedge b$ and $e:=a\vee b$.
By Lemma~\ref{lem:component-congruences},
\[
\negt d\leqslant\negt e,
\qquad
\nego d\leqslant\nego e.
\]
By Lemma~\ref{lem:lindenbaum-structure}, the operation
$k:=\negt\circ\nego$ is antitone and the negations are commuting involutions.
Applying $k$ to the first inequality gives
\[
\nego e=k\negt e\leqslant k\negt d=\nego d.
\]
Together with the second inequality, this yields $\nego d=\nego e$.
Involutivity gives $d=e$.
Since $d\leqslant a,b\leqslant e$, we obtain $a=b$.
Thus $\theta_t\cap\theta_1\subseteq\Delta_L$.
The reverse inclusion follows from reflexivity of both congruences.

The operations on $L_t\times L_1$ are defined coordinatewise.
Since both $q_t$ and $q_1$ are homomorphisms, so is $q$.
Moreover,
\[
q(a)=q(b)
\iff (a\theta_t b\ \text{and}\ a\theta_1 b)
\iff a=b,
\]
so $q$ is injective.
For each coordinate projection $\pi_i$, surjectivity of $q_i$ gives
$\pi_i(q[L])=q_i[L]=L_i$. Hence $q$ is a subdirect embedding.

Finally, suppose $q(a)\leqslant q(b)$.
Then $q(a\wedge b)=q(a)\wedge q(b)=q(a)$.
Injectivity gives $a\wedge b=a$, hence $a\leqslant b$.
Thus $q$ reflects order.
\end{proof}

\subsection{Decomposition and completeness}\label{sec:decomposition}

\begin{theorem}\label{thm:syntactic-decomposition}
The derivability relations satisfy
\[
\dern\ =\ \dert\cap\derone.
\]
\end{theorem}

\begin{proof}
The inclusion $\dern\subseteq\dert\cap\derone$ follows from
Corollary~\ref{cor:embedding}.
Conversely, suppose $A\dert B$ and $A\derone B$.
By Lemma~\ref{lem:component-quotients},
\[
q_t([A])\leqslant q_t([B]),
\qquad
q_1([A])\leqslant q_1([B]).
\]
The order on $L_t\times L_1$ is defined coordinatewise, so
$q([A])\leqslant q([B])$.
By Lemma~\ref{lem:component-separation}, $q$ reflects order,
hence $[A]\leqslant[B]$.
By the definition of the order on $L$, this means $A\dern B$.
\end{proof}

\begin{corollary}\label{cor:completeness-n}
For all $A,B\in\For$,
\[
A\dern B\iff A\entn B.
\]
\end{corollary}

\begin{proof}
Soundness is Theorem~\ref{thm:soundness-n}.
Conversely, if $A\entn B$, then $A\entt B$ and $A\entone B$
by Definition~\ref{def:entailment}.
Theorem~\ref{thm:completeness-t} and Remark~\ref{rem:mirror} yield
$A\dert B$ and $A\derone B$, so $A\dern B$ by
Theorem~\ref{thm:syntactic-decomposition}.
\end{proof}

\section{Finite models and decidability}\label{sec:finite-models}

\begin{lemma}\label{lem:finite-order-completion}
Let $S$ be a nonempty finite partially ordered set with an antitone
involution $n$. There is a finite complete lattice $K$ with an antitone
involution \({}^\ast\) and an order embedding $e\colon S\to K$.
The embedding preserves all meets and joins that exist in $S$ and
satisfies $e(n(y))=e(y)^\ast$ for every $y\in S$.
\end{lemma}

\begin{proof}
For $T\subseteq S$, write $T^\ell$ and $T^u$ for the sets of common
lower and upper bounds of $T$ in $S$. Let
$K=\mathrm{DM}(S)=\{T^\ell\mid T\subseteq S\}$ be the
Dedekind--MacNeille completion, ordered by inclusion. It is a complete
lattice and is finite because its elements are subsets of $S$.
The canonical map $e(y)=\{y\}^\ell$ is an order embedding that preserves
all meets and joins that exist in $S$ \cite[Theorem~7.40]{DaveyPriestley}.
Indeed, if $m$ is the meet of $T\subseteq S$, then the intersection of
the principal lower sets $e(t)$ for $t\in T$ is $T^\ell=e(m)$.
For joins, if $j$ is the join of $T$ in $S$, every upper bound
$D=U^\ell$ of the sets $e(t)$ contains $T$. Thus each member of $U$
is an upper bound of $T$ and lies above $j$. It follows that $e(j)$
is contained in $D$. Since $e(j)$ itself contains all the $e(t)$,
it is their join in $K$.

For $D\in K$, put $D^\ast=(n[D])^\ell$, where $n[D]$ is the image of
$D$ under $n$. This defines an antitone map on $K$, since taking common
lower bounds reverses inclusion. As $n$ is an antitone involution,
$x\leqslant n(d)$ if and only if $d\leqslant n(x)$. Thus $x\in D^\ast$
if and only if $n(x)\in D^u$, so $n[D^\ast]=D^u$.

For $D=T^\ell\in K$, every element of $T$ is an upper bound of $D$.
Hence $(D^u)^\ell\subseteq T^\ell=D$. The reverse inclusion holds
because each element of $D$ lies below every upper bound of $D$.
It follows that $(D^\ast)^\ast=(D^u)^\ell=D$.

Since $n$ reverses and reflects the order, $n[e(y)]=\{n(y)\}^u$.
This set has least element $n(y)$, so its common lower bounds form
$e(n(y))$. Hence $e(y)^\ast=e(n(y))$.
\end{proof}

\begin{lemma}\label{lem:finite-frame-realisation}
Let $K$ be a finite nontrivial lattice with an antitone involution
\({}^\ast\). There is a finite $\NDst$-frame with
$\Lstab_t=\Lstab_1$ and a lattice embedding
$\beta\colon K\to\Lstab_t$ satisfying
\[
\beta(a^\ast)=\Nc{i}(\beta(a))
\qquad(a\in K,\ i\in\{t,1\}).
\]
\end{lemma}

\begin{proof}
Order the disjoint filter--ideal pairs of $K$ by componentwise inclusion.
They form a finite poset containing $(\{1_K\},\{0_K\})$.
Its set $X_K$ of maximal elements is therefore finite and nonempty.
For $x=(F,I)$ and $y=(G,J)$ in $X_K$, put $x\preceq_F y$ if
$F\subseteq G$ and $x\preceq_I y$ if $I\subseteq J$.

For $U\subseteq K$, write $U^\ast=\{u^\ast\mid u\in U\}$.
On the set of all disjoint filter--ideal pairs, define
$N(F,I)=(I^\ast,F^\ast)$. The antitone involution of $K$ interchanges
filters and ideals and preserves disjointness by injectivity.
The map $N$ preserves componentwise inclusion and is an involution
on these pairs.
If $(F,I)\in X_K$, a strict extension of $N(F,I)$ as a disjoint pair
would give a strict extension of $(F,I)$ on applying $N$ again.
Thus $N$ maps $X_K$ to itself and is an involution there.
Moreover, $F\subseteq G$ if and only if $F^\ast\subseteq G^\ast$, so
$x\preceq_F y$ if and only if $N(x)\preceq_I N(y)$.

Use $(\preceq_F,\preceq_I)$ as the pair of orders in each component
and set $\Nt=\No=N$. The two stars commute because they coincide.
By Definition~\ref{def:frame}, this gives a finite $\NDst$-frame.
The two lattices of stable sets coincide because the pairs of orders
are the same.

Urquhart's representation in \eqref{eq:urqrep} gives a lattice embedding
into this common lattice by
\[
\beta(a)=\{(F,I)\in X_K\mid a\in F\}.
\]
By \eqref{eq:urqsecond},
$\rmap{i}(\beta(a))=\{(F,I)\in X_K\mid a\in I\}$ for each $i\in\{t,1\}$.
For any $x=(F,I)\in X_K$, we have
\[
\begin{aligned}
x\in N[\rmap{i}(\beta(a))]
&\iff N(x)\in\rmap{i}(\beta(a))\\
&\iff a\in F^\ast
 \iff a^\ast\in F
 \iff x\in\beta(a^\ast).
\end{aligned}
\]
The first equivalence uses $N^2=\mathrm{id}_{X_K}$. Since
$\Nc{i}(\beta(a))=N[\rmap{i}(\beta(a))]$, this proves the required identity.
\end{proof}

\begin{theorem}\label{thm:finite-model-property}
Each of $\entt$, $\entone$ and $\entn$ is determined by the models on
finite $\NDst$-frames.
\end{theorem}

\begin{proof}
Fix $i\in\{t,1\}$ and suppose $A\not\models^{i}_{oe}B$. Choose a model
in which $\Pset{A}{i}\nsubseteq\Pset{B}{i}$. Let $\Sigma$ be the set of
subformulae of $A$ and $B$. In the lattice $\Lstab_i$ of stable sets, put
\[
S_0=\{\Pset{C}{i}\mid C\in\Sigma\},
\qquad
S=S_0\cup\Nc{i}[S_0].
\]
The set $S$ is finite and contains $\Pset{A}{i}$ and $\Pset{B}{i}$.
Order $S$ by inclusion. The antitone involution $\Nc{i}$ of $\Lstab_i$
maps $S$ onto itself by construction. Apply
Lemma~\ref{lem:finite-order-completion} to $S$ with the restriction of
$\Nc{i}$. This gives a finite complete lattice $K$ with an antitone
involution \({}^\ast\) and an order embedding $e\colon S\to K$.

For each variable $p$ occurring in $A$ or $B$, set $v(p)=e(\Pset{p}{i})$.
Assign the bottom element of $K$ to the other variables. Extend $v$ to
formulae using the meet and join of $K$, with $\neg_i$ interpreted by
\({}^\ast\) and the other negation by the identity.

By Corollary~\ref{cor:homomorphism}, the truth set of a conjunction or
disjunction in $\Sigma$ is the meet or join in $\Lstab_i$ of the truth sets
of its immediate subformulae. All three truth sets belong to $S$.
In the case of a conjunction, any lower bound in $S$ is also a lower
bound in $\Lstab_i$ and is therefore below the truth set of the
conjunction. The dual argument applies to a disjunction. Thus these
meets and joins exist already in $S$ and are preserved by $e$.

We now show $v(C)=e(\Pset{C}{i})$ by induction on $C\in\Sigma$.
For variables this is the definition of $v$. For binary connectives
it follows from preservation of the indicated meets and joins.
For $C=\neg_i D$, the induction hypothesis gives
$v(D)=e(\Pset{D}{i})$. Preservation of the involution then gives
$v(C)=e(\Nc{i}(\Pset{D}{i}))=e(\Pset{C}{i})$.
The other negation acts as the identity in component $i$ and in $K$.
Since $e$ reflects order,
\[
v(A)=e(\Pset{A}{i})\nleqslant e(\Pset{B}{i})=v(B).
\]
In particular, $K$ is nontrivial.

Apply Lemma~\ref{lem:finite-frame-realisation} to $K$. With $i$ fixed as
at the start of the proof, let $j$ be the other component index.
On the finite frame, assign
to each variable $p$ the truth sets
\[
\Pset{p}{i}=\beta(v(p)),
\qquad
\Pset{p}{j}=\beta(0_K).
\]
For $k\in\{i,j\}$, the corresponding falsity set is $\rmap{k}(\Pset{p}{k})$.
Since $\Pset{p}{k}$ is $\lmap{k}$-stable, the pair
$(\Pset{p}{k},\rmap{k}(\Pset{p}{k}))$ is $k$-conjugate by
Lemma~\ref{lem:conjugate}. These pairs define a valuation.
The choice in component $j$ supplies the second conjugate pair for each
variable.

By Corollary~\ref{cor:homomorphism}, evaluation in component $i$ uses
only the values assigned in that component. Conjunction and
disjunction are interpreted by meet and join, $\neg_i$ by $\Nc{i}$, and
$\neg_j$ by the identity. The embedding $\beta$ preserves meets and joins
and satisfies $\beta(a^\ast)=\Nc{i}(\beta(a))$. Thus induction on $C$ gives
$\Pset{C}{i}=\beta(v(C))$ in the new model.

Since $v(A)\nleqslant v(B)$ and $\beta$ reflects order,
$\beta(v(A))\nsubseteq\beta(v(B))$. A point in
$\beta(v(A))\setminus\beta(v(B))$ therefore witnesses the failure
of $A\models^i_{oe}B$ on the finite frame.

If $A\not\entn B$, then $A\not\models^{i}_{oe}B$ for at least one
$i\in\{t,1\}$ by Definition~\ref{def:entailment}. The finite countermodel
for that component also witnesses failure of the joint entailment.
\end{proof}

\begin{corollary}\label{cor:decidability}
The three entailment relations and their corresponding derivability
relations are decidable.
\end{corollary}

\begin{proof}
Fix formulae $A,B$ and one of the three entailment relations. Write
$\vdash$ for derivability in its calculus. Using an effective numbering
of the propositional variables and the axiom and rule schemata,
enumerate all finite derivations in that calculus. Search this
enumeration for a derivation of $A\vdash B$.

To search for a finite countermodel, take $X=\{1,\ldots,m\}$ for each
$m\geqslant 1$ and test the conditions of Definition~\ref{def:frame} on
all choices of four binary relations and two maps $X\to X$.
Every finite frame can be relabelled onto one of these carriers.
On each resulting frame, enumerate the assignments of a $t$-conjugate
pair and a $1$-conjugate pair to each variable occurring in $A$ or $B$,
testing conjugacy by \eqref{eq:conjugate}. Each assignment extends to
a valuation by assigning $(\varnothing,X)$ in both components to the
remaining variables. In each component $k$, $\rmap{k}(\varnothing)=X$,
and reflexivity gives $\lmap{k}(X)=\varnothing$.

The truth and falsity sets of $A$ and $B$ depend only on the values
of their variables, by structural induction on the semantic conditions in
Definition~\ref{def:forc}. On a finite carrier those conditions give a
recursive computation of the two pairs of sets for each formula.
Testing the chosen entailment in this model therefore requires only
one or two inclusions, according to Definition~\ref{def:entailment}.
There are finitely many frames and assignments to test for each $m$,
so the search reaches every finite countermodel up to relabelling and
agreement on the variables in $A$ and $B$.

Alternate steps of the two searches. By Theorem~\ref{thm:completeness-t},
Remark~\ref{rem:mirror} and Corollary~\ref{cor:completeness-n}, $A\vdash B$
if and only if the chosen entailment holds. If it holds, the first search
finds a derivation of $A\vdash B$. Otherwise,
Theorem~\ref{thm:finite-model-property} ensures that the second search
finds a finite countermodel. The search that terminates first determines
whether $A\vdash B$ and whether the corresponding entailment holds.
\end{proof}

\section{Conclusion}\label{sec:conclusion}

The four- and nine-valued entailments arise on specific frames within a
class that also admits nondistributive lattices of stable sets. For the
full class, the subdirect representation of the joint Lindenbaum
algebra gives the syntactic decomposition. Joint completeness follows
from this decomposition and completeness of the two component calculi.
To obtain a finite countermodel, we complete the ordered set of
subformula values and their involutive images in a component where
entailment fails. The completion preserves the involution and the
meets and joins used to evaluate the subformulae. We then represent
the resulting lattice on a finite frame. If joint entailment fails,
it fails in at least one component. A finite countermodel in that
component therefore suffices.
Together with the effective calculi, the finite model property gives
decidability for all three entailments.

Corollary~\ref{cor:finiteentailment} identifies the entailments on
$\Frm_2$ and $\Frm_3$ with the four- and nine-valued entailments,
respectively. A further task is to axiomatise each of these by adding
postulates to the joint calculus.

By Lemma~\ref{lem:reduction}, in each component forward preservation
of truth holds at all points of a model if and only if backward
preservation of falsity does.
In $\mathrm{SIXTEEN}_3$, Shramko and Wansing~\cite{ShramkoWansing16}
define separate entailments using the truth and falsity orders.
Odintsov~\cite[Section~4]{OdintsovSixteen} presents each as an
intersection of coordinate entailments. Odintsov and
Wansing~\cite[Theorem~4.1]{OdintsovWansingBilattices} give a calculus for
both in the implication-free language using commutative distributive
bilattices with conflation. A separate task for our semantics is to
weaken conjugacy so that the two preservation conditions yield distinct
entailments within each component while keeping each semi-negation
transparent in the other component.

\section*{Acknowledgments}

The author used OpenAI Codex to locate sources and verify bibliographic
details, prepare and edit portions of the English text, and format the
manuscript in LaTeX.
The author takes full responsibility for the proofs and the final text of
the manuscript.

\end{document}